\UseRawInputEncoding
\RequirePackage{fix-cm}
\documentclass[smallcondensed]{svjour3}     
\usepackage{amsmath}
\usepackage[numbers,sort&compress]{natbib}
\usepackage{amsfonts,amssymb,amsmath}
\usepackage{mathrsfs}
\usepackage{color}
\usepackage[colorlinks,bookmarksopen,bookmarksnumbered,citecolor=black, linkcolor=blue, urlcolor=black]{hyperref}
\newcommand{\beq}{\begin{equation}}
\newcommand{\enq}{\end{equation}}
\newtheorem{Theorem}{Theorem}
\newtheorem{Lemma}{Lemma}
\newtheorem{Corollary}{Corollary}
\newtheorem{Definition}{Definition}
\newtheorem{Remark}{Remark}

\newcommand{\benu}{\begin{enumerate}}
\newcommand{\beqa}{\begin{eqnarray}}
\newcommand{\beqan}{\begin{eqnarray*}}
\newcommand{\eay}{\end{array}}
\newcommand{\edm}{\end{displaymath}}
\newcommand{\eenu}{\end{enumerate}}
\newcommand{\eeq}{\end{equation}}
\newcommand{\eeqa}{\end{eqnarray}}
\newcommand{\eeqan}{\end{eqnarray*}}

\newcommand{\br}{\begin{Remark}}
\newcommand{\er}{\end{Remark}}

\newcommand{\bqa}{\begin{eqnarray}}
\newcommand{\eqa}{\end{eqnarray}}
\newcommand{\bqw}{\begin{eqnarray*}}
\newcommand{\eqw}{\end{eqnarray*}}

\newcommand{\bea}{\begin{array}{cc}}
\newcommand{\ena}{\end{array}}

\begin{document}

\title{Existence and regularity of pullback attractors for a non-autonomous diffusion equation with delay and nonlocal diffusion in time-dependent spaces}


\author{Yuming Qin         \and
        Bin Yang 
}

\institute{Yuming Qin \at
              Department of  Mathematics, Institute for Nonlinear Science, Donghua University, Shanghai, 201620, P. R. China. \\
              Tel.: +0086-21-67874272\\
              \email{yuming$\_$qin@hotmail.com}           
           \and
           Bin Yang \at
              College of Information Science and Technology, Donghua University, Shanghai, 201620, P. R. China.\\
               \email{binyangdhu@163.com}
}

\date{Received: date / Accepted: date}
\maketitle

\large{
\begin{abstract}
In this paper, we study the asymptotic behavior of solution to a non-autonomous diffusion equations with delay containing some hereditary characteristics and nonlocal diffusion in time-dependent space $C_{\mathcal{H}_{t}(\Omega)}$. When the nonlinear function $f$ satisfies the polynomial growth of arbitrary order $p-1$ $(p \ge 2)$ and the external force $h \in L_{l o c}^{2}\left(\mathbb{R} ; H^{-1}(\Omega)\right)$, we establish the existence and regularity of the time-dependent pullback attractors.
\end{abstract}

\section{\large Introduction}
In recent years, a lot of scholars have been keen to discuss the asymptotic behavior of solutions to the diffusion equations. The assumption for function $a(\cdot)$ of problem $u_{t}-a(l(u))\Delta u=f$ in Lovat \cite{l.4} is only $0<\bar m \leq a(s) \leq \bar M$ for any positive constants $\bar m$ and $\bar M$, which avoids that weak solutions of the diffusion equations only exist in a finite-time interval. In addition, Chipot and Zheng \cite{cz.4} showed that due to the existence of the nonlocal function $a(\cdot)$, it is impossible to guarantee the existence of a Lyapunov structure. Since diffusion equations with nonlocal diffusion are widely used in ecology, epidemiology, materials science, neural network and other disciplines, and scientific researches toward them are challenging and forward-looking, they have been extensively studied (see \cite{chm2,chm16,qy.2,y.4}).

In this paper, we consider the following non-autonomous diffusion equation with delay and nonlocal diffusion
\begin{equation}
\left\{\begin{array}{ll}
\partial_{t}u-\varepsilon(t) \partial_{t}\Delta u-a(l(u)) \Delta u=f(u)+g(t, u_{t})+h(t) & \text { in } \Omega \times(\tau, \infty), \\
u(x,t)=0 & \text { on } \partial \Omega\times(\tau, \infty), \\
u(x, \tau+\theta)=\phi(x, \theta),  &\,\, x \in \Omega,\, \theta \in[-k, 0],
\end{array}\right.\label{1.1-4}
\end{equation}
in time-dependent space $C_{\mathcal{H}_{t}(\Omega)}$, where $\Omega  \subset \mathbb{R}^{n}\,(n \ge 3)$ is a bounded domain with smooth boundary $\partial \Omega$, the initial time $\tau  \le t \in \mathbb R$, $\phi \in C\left([-k, 0] ; \mathcal{H}_{t}(\Omega)\right)$ is the initial datum, $k\,(>0)$ is the length of the delay effects, $f$ is a nonlinear function, $g$ is a delay operator containing some hereditary characteristics and we suppose $u_{t}$ is defined in $[-k, 0]$ and satisfies $u_{t}(\theta)=u(t+\theta)$. Meanwhile, let the external force $h(x, t) \in L_{\mathrm{loc}}^{2}(\mathbb R ; H^{-1}(\Omega))$.

Suppose $C([-k, 0] ; X)$ is a Banach space equipped with the sup-norm $\|\cdot\|_{C([-k,0];X)}$ and metric $d_{C}(\cdot, \cdot)$, and denote it as $C_{X}$. The time-dependent space $\mathcal H_{t}(\Omega)$ is equipped with norm $\|u\|_{\mathcal{H}_{t}(\Omega)}^{2}=\|u\|_{2}^{2}+\varepsilon(t)\|\nabla u\|_{2}^{2}$. For brevity, the norm of $L^{2}(\Omega)$ in this paper is denoted as $\|\cdot\|_{2}=\|\cdot\|$. Besides, we assume $C_{L^{2}(\Omega)}$ is a Banach space with sup-norm, that is, for any $u \in C_{L^{2}(\Omega)}$, its norm is defined as $\|u\|_{C_{L^{2}(\Omega)}}=\max\limits _{t \in[-k, 0]}\|u\|_{2}$. Similarly, for any $u \in C_{H_{0}^{1}(\Omega)}$, we define its norm as $\|u\|_{C_{H_{0}^{1}(\Omega)}}=\max\limits _{t \in[-k, 0]}(\|u\|_{2}+\|\nabla u\|_{2})$. For any $t \in \mathbb{R}$, the time-dependent space $C_{\mathcal{H}_{t}(\Omega)}$ is endowed with the norm
$$
\|u\|_{C_{\mathcal{H}_{t}(\Omega)}}^{2}=\|u\|_{C_{L^{2}(\Omega)}}^{2}+| \varepsilon_{t}| \|\nabla u\|_{C_{L^{2}(\Omega)}}^{2},
$$
where $\varepsilon_{t}=\varepsilon(t+\theta)$ with $\theta \in[-k, 0]$ and $| \varepsilon_{t}|$ is the absolute value of $\varepsilon_{t}$.

Assume $\varepsilon(t) \in C^{1}(\mathbb{R})$ is a time-dependent function satisfying
\begin{equation}
\lim _{t \rightarrow+\infty} \varepsilon(t)=1,
\label{1.2-4}
\end{equation}
and there exists a constant $L>0$ such that
\begin{equation}
\sup _{t \in \mathbb{R}}(|\varepsilon(t)|+|\varepsilon^{\prime}(t)|) \leq L.
\label{1.3-4}
\end{equation}

In addition, the function $a(l(u)) \in C\left(\mathbb{R}; \mathbb{R}^{+}\right)$ is the nonlocal diffusion term of equation $(\ref{1.1-4})$ and satisfies
\begin{equation}
\frac{1}{2}\left(3+L+\varepsilon^{\prime}(t)\right)<m \leq a(s) \leq M,  \quad \forall \, s \in \mathbb{R},
\label{1.4-4}
\end{equation}
where $m$ and $M$ are positive constants.
Assume $l(u): L^{2}(\Omega) \to \mathbb R$ is a continuous linear functional acting on $u$ that satisfies for some $j \in L^{2}(\Omega)$,
\begin{equation}
l(u)=l_{j}(u)=\int_{\Omega} j(x) u(x) d x.
\label{1.5-4}
\end{equation}

Besides, suppose the nonlinear function $f \in C^{1}(\mathbb{R}, \mathbb{R})$ and satisfies
\begin{equation}
(f(u)-f(v))(u-v) \leq \tilde{\eta}(u-v)^{2}, \quad \forall \, u, v\in \mathbb R,
\label{1.6-4}
\end{equation}
and
\begin{equation}
-C_{0}-C_{1}|u|^{p} \leq f(u) u \leq C_{0}-C_{2}|u|^{p}, \quad p \geq 2,
\label{1.7-4}
\end{equation}
for some positive constants $\tilde{\eta}$, $C_{0}$, $C_{1}$ and $C_{2}$. Let $\mathcal F(u)=\int_{0}^{u} f(r) d r$, the same as in \cite{gm.4}, then from $(\ref{1.6-4})-(\ref{1.7-4})$, it follows that there exists some positive constants ${\widetilde{C}_{i}}$ $(i=1, 2)$ such that
\begin{equation}
-{\widetilde{C}_{0}}-{\widetilde{C}_{1}}|u|^{p} \leq \mathcal{F}(u) \leq {\widetilde{C}_{0}}-{\widetilde{C}_{2}}|u|^{p}.
\label{1.8-4}
\end{equation}

Furthermore, assume the delay operator $g: \mathbb{R} \times C_{L^{2}(\Omega)} \rightarrow L^{2}(\Omega)$ and it satisfies

(A1) for any $\nu \in C_{L^{2}(\Omega)}$, the function $\mathbb{R} \ni t \mapsto g(t, \nu) \in L^{2}(\Omega)$ is measurable;

(A2) $g(t, 0)=0,$ for all $t \in \mathbb{R}$;

(A3) there exists $C_{g}>0$ such that for all $t \in \mathbb{R}$ and $\nu_{1}, \nu_{2} \in C_{L^{2}(\Omega)}$, it holds
$$
\|g(t, \nu_{1})-g(t, \nu_{2})\|^{2} \leq C_{g}\|\nu_{1}-\nu_{2}\|^{2}_{C_{L^{2}(\Omega)}}.
$$

Throughout this paper, the inner product of $L^{2}(\Omega)$ is denoted as $(\cdot, \cdot)$ and the norms of $L^{\gamma}(\Omega)$ with $\gamma \in \mathbb R$, $H_{0}^{1}(\Omega)$ and $H^{-1}(\Omega)$ are denoted by ${\|\cdot\|}_{\gamma}$, ${\|\cdot\|}_{1}$ and ${\|\cdot\|}_{-1}$, respectively.

There are many studies for the diffusion equations with delays related to problem $(\ref{1.1-4})$. Hu and Wang \cite{hw.4} investigated the existence of pullback attractors for $\partial_{t} u-\Delta \partial_{t} u-\Delta u=f(t, u(t-\rho(t)))+g(t)$  in $C_{H_{0}^{1}(\Omega)}$ and $C_{H^{2}(\Omega) \cap H_{0}^{1}(\Omega)}$, where $\rho$ is a delay function and $f$ contains some memory effects in a fixed time interval with length $h>0$. Later on, Caraballo and M\'{a}rquez-Dur\'{a}n \cite{cm.4} proved the existence and eventual uniqueness of stationary solutions, as well as their exponential stability for $\partial_{t} u-\Delta\partial_{t} u-\Delta u=g\left(t, u_{t}\right)$ in $L^{2}(\Omega)$. Besides,  Garc\'{i}a-Luengo and Mar\'{i}n-Rubio \cite{gm.4} obtained two different families of the minimal pullback attractors for $\partial_{t} u-\Delta u=f(u)+g\left(t, u_{t}\right)+k(t)$ in $H^{-1}(\Omega)$. Furthermore, Zhu and Sun \cite{zs2} verified the existence of pullback attractors for this equation in $C_{H_{0}^{1}(\Omega)}$.

In addition, Wang and Kloeden \cite{wk.4} studied the existence of a uniform attractor in $L^{2}(\Omega)$ for the multi-valued process associated with ${\partial_{t} u}-\Delta u+\lambda u=f\left(x, u_{t}\right)+g(t, x)$. Caraballo, M\'{a}rquez-Dur\'{a}n and Rivero \cite{cmr.4} established the existence of the pullback attractors for ${\partial_{t} u}-\gamma(t) \Delta {\partial_{t} u}-\Delta u=g(u)+f\left(t, u_{t}\right)$ in $C_{H_{0}^{1}(\Omega)}$, where $\gamma: \mathbb{R} \rightarrow$ $(0,+\infty)$ is a continuous bounded function with $0<\gamma_{0} \leq \gamma(t) \leq \gamma_{1}<\infty$. Moreover, Harraga and Yebdri \cite{hy.4} proved the existence of the pullback $\mathcal D$-attractors in $H_{0}^{1}(\Omega)$ for ${\partial_{t}}u-\Delta {\partial_{t}}u-\Delta u=b(t, u(t-\rho(t))(x)+g(t, x)$. Zhu, Xie and Zhang \cite{zxz.4} studied the asymptotic behavior of pullback attractors for $\partial_{t} u-\varepsilon(t) \partial_{t} \Delta u-\Delta u+f(u)=g\left(t, u_{t}\right)+k(x)$ in time-dependent space $\mathcal H_{t}(\Omega)$. Additionally, some scholars have also considered the well-posedness of solutions to the reaction-diffusion equations with other forms of delays (see \cite{cmr2.4,lt.4,sc.4}).

On the basis of Garc\'{i}a-Luengo and Mar\'{i}n-Rubio \cite{gm.4}, we follow some of the methods in Zhu and Sun \cite{zs2} to consider the existence of the time-dependent pullback $\mathcal D_{{\eta}_{1}}$-attractor of problem $(\ref{1.1-4})$ in time-dependent space $C_{\mathcal{H}_{t}(\Omega)}$. When compared to their equations, the terms $a(l(u))$ and $\varepsilon(t)$ bring some difficulties to verify the desired results, and now we will illustrate them and explain our strategies and innovations.

$(1)$ In \cite{y.4}, we consider the existence and upper semicontinuity of problem (\ref{1.1-4}) without the term $g(t,u_{t})$. Moreover, many scholars have studied the attractors in time-dependent space $\mathcal H_{t}(\Omega)$ (see \cite{chm,chm4.2,mwx,qy.2}). To control the term $-\varepsilon^{\prime}(t)$ during energy estimations, they basically required the time-dependent term $\varepsilon(t)$ to be a decreasing function and using the transform
$$
\varepsilon(t) \frac{d}{d t}\|\nabla u\|^{2}=\frac{d}{d t}\left(\varepsilon(t)\|\nabla u\|^{2}\right)-\varepsilon^{\prime}(t)\|\nabla u\|^{2}
$$
to deal with the term $-\varepsilon^{\prime}(t)\|\nabla u\|^{2}$, then they can easily derive the range of weak solutions in $\mathcal H_{t}(\Omega)$. In this paper, we restrict the lower bound of the nonlocal function $a(\cdot)$ without supposing the increase or decrease of $\varepsilon(t)$, which is a brand new attempt and weaker than the conditions in \cite{chm,chm4.2,mwx,qy.2}.

$(2)$ The phase spaces of \cite{gm.4} and \cite{zs2} are $H^{-1}(\Omega)$ and $C_{H_{0}^{1}(\Omega)}$, respectively, while our phase space $C_{\mathcal{H}_{t}(\Omega)}$ is a time-dependent space, and since problem $(\ref{1.1-4})$ contains the nonlocal function $a(l(u))$ and the time-dependent function $\varepsilon(t)$, they ensure that problem $(\ref{1.1-4})$ to have wider applications in real life. Meanwhile, these terms make it impossible to directly conclude the specific structure of absorbing set and the asymptotic compactness of the process by general energy estimations, which requires a series of delicate and multiple calculations. To this end, as in \cite{zs2}, we shall use a special Gronwall lemma (see \cite{Evans,Q1,Q2,Q3}) to conclude our desired results, which will be shown in Lemma \ref{lem4.1-4}.

$(3)$ In \cite{gm.4} and \cite{y.4}, it is difficult to establish regularity for process due to the term $-\partial_{t} \Delta u$. In problem $(\ref{1.1-4})$, this term is replaced by $-\varepsilon(t)\partial_{t} \Delta u$, and we also obtain the regularity of the process in $C_{\mathcal H_{t}(\Omega)}$ by decomposition method and using some calculations and estimates similar to those used to demonstrate the existence of the time-dependent pullback $\mathcal D_{\eta_{1}}$-attractor, which provides some ideas for studying pullback attractors in spaces of higher regularity.

This paper is organized as follows. In $\S 2$, we shall introduce some useful abstract definitions, theorems and lemmas. Next, we shall prove the existence and uniqueness of weak solutions to problem $(\ref{1.1-4})$ by Faedo-Galerkin approximations in $\S 3$. Furthermore, in $\S 4$ we shall verify the existence of the time-dependent pullback $\mathcal D_{{\eta}_{1}}$-attractor ${\mathcal{A}}_{{\eta}_{1}}$ in time-dependent space $C_{\mathcal{H}_{t}(\Omega)}$. Finally, we shall derive the regularity of ${\mathcal{A}}_{{\eta}_{1}}$ in $\S 5$.
\section{\large Preliminaries}
In this section, we will introduce some needed abstract concepts, such as definitions and properties of function spaces and attractors.

The time-dependent space $C_{\mathcal{H}_{t}^{1}(\Omega)}$, more regular than $C_{\mathcal{H}_{t}(\Omega)}$, is endowed with the norm
$$
\|u\|_{C_{\mathcal{H}_{t}^{1}(\Omega)}}^{2}=\|\nabla u\|_{C_{L^{2}(\Omega)}}^{2}+| \varepsilon_{t}| \|\Delta u\|_{C_{L^{2}(\Omega)}}^{2},
$$
where $\varepsilon_{t}=\varepsilon(t+\theta)$ with $\theta \in[-k, 0]$.

The closed $R$-ball with the origin as the center and $R$ as the radius in $C_{X}$ is denoted as $${\mathcal{\bar B}_{{C_{X}}}}(0, R) = \left\{ {u \in {C_{X}}:\left\| u \right\|_{{C_{X}}}^2 \le R} \right\}.$$

The Hausdorff semidistance between two nonempty  sets $A_{1}, A_{2} \subset C_{X}$ is denoted by
$$
dist_{C_{X}}(A_{1},A_{2})=\sup _{x \in A_{1}} \inf _{y \in A_{2}}\|x-y\|_{C_{X}} \, .
$$

\begin{Definition} {\rm(\cite{mwx,zxz})}
Assume $\left\{X_{t}\right\}_{t \in \mathbb{R}}$ is a family of normed spaces. A process or a two-parameter semigroup on $\left\{X_{t}\right\}_{t \in \mathbb{R}}$ is a family $\{U(t, \tau) \mid t, \tau \in \mathbb{R}, t \ge \tau\}$ of mapping $U(t, \tau): X_{\tau} \rightarrow X_{t}$ satisfies that $U(\tau, \tau)u=u$ for any $u \in X_{\tau}$ and $U(t, s) U(s, \tau)=U(t, \tau)$ for all $t \geq s \ge \tau$.
\label{def2.1-4}
\end{Definition}

\begin{Definition} {\rm(\cite{chm})}
For any $\sigma>0$, let $\mathcal D$ be a nonempty class of all families of parameterized sets $\widehat{D}=\left\{D(t): t \in \mathbb{R}\right\} \subset \Gamma(X_{t})$ such that
$$
\lim _{\tau \rightarrow-\infty}\left(e^{\sigma \tau} \sup _{u \in D(\tau)}\|u\|_{X_{t}}^{2}\right)=0,
$$
where $\Gamma(X_{t})$ denotes the family of all nonempty subsets of $\left\{X_{t}\right\}_{t \in \mathbb{R}},$ then $\mathcal D$ will be called a tempered universe in $\Gamma(X_{t})$.
\label{def2.2-4}
\end{Definition}

\begin{Definition} {\rm(\cite{psz,zs2})}
A process ${\{ U(t,\tau )\} _{t \ge \tau }}$ is pullback $\mathcal D$-asymptotically compact on $\left\{X_{t}\right\}_{t \in \mathbb{R}}$, if for any $t \in \mathbb{R}$, any ${\widehat D} \in {{\cal D}}$, any sequence ${\left\{ {{\tau _n}} \right\}_{n\in{\mathbb{N}^ + }}}\subset( - \infty ,t]$ and any sequence ${\left\{ {{x_n}} \right\}_{n \in {\mathbb{N}^ + }}} \subset {D({\tau}_{n})} \subset {{ X}_t}$, the sequence $\left.\left\{U (t, \tau\right) x_{n}\right\}_{n\in{\mathbb{N}^ + }}$ is relatively compact in $X_{t}$ when ${\tau _n} \to  - \infty $.
\label{def2.3-4}
\end{Definition}

\begin{Definition} {\rm(\cite{psz,zs2})}
A family ${\widehat D_{0}=\left\{D_{0}(t): t \in \mathbb{R}\right\}} \subset \Gamma(X_{t})$ is pullback $\mathcal D$-absorbing for the process ${\{ U(t,\tau )\} _{t \ge \tau }}$ on $\left\{X_{t}\right\}_{t \in \mathbb{R}}$, if for any $t \in \mathbb{R}$ and ${\widehat D} \in {{\cal D}}$, there exists a ${\tau _0} = {\tau _0}(t,{\widehat D}) < t$ such that
$
U(t, \tau) D(\tau) \subset D_{0}(t)
$
for any $\tau \leq \tau_{0}(t, \widehat{D})$.
\label{def2.4-4}
\end{Definition}

\begin{Definition} {\rm(\cite{psz,zs2})}
A family ${\mathcal{A}_{\eta}}=\{\mathcal{A}_{\eta}(t): t \in \mathbb{R}\} \subset \Gamma\left(X_{t}\right)$
is called a time-dependent pullback $\mathcal D$-attractor for the process ${\{ U(t,\tau )\} _{t \ge \tau }}$ on $\left\{X_{t}\right\}_{t \in \mathbb{R}}$, if the following properties hold:

(i) the set $\mathcal{A}_{\eta}(t)$ is compact in $X_{t}$ for any $t \in \mathbb{R}$;

(ii) ${\mathcal{A}}_{\eta}$ is pullback $\mathcal D$-attracting in $X_{t}$, i.e.,
$$
\lim _{\tau \rightarrow-\infty} {dist}_{X_{t}}\left(U(t, \tau) D(\tau), \mathcal A_{\eta}(t)\right)=0,
$$
for any ${\widehat D} \in {{\cal D}}$ and $t \in \mathbb{R}$;

(iii) ${\mathcal{A}}_{\eta}$ is invariant, i.e., $U(t, \tau)\mathcal{A}_{\eta}(\tau)={\mathcal{A}_{\eta}(t)}$, for any $\tau \leq t$.
\label{def2.5-4}
\end{Definition}

Furthermore, we will introduce the following definitions and lemmas, which contribute to prove the existence of the time-dependent pullback $\mathcal D$-attractors.

\begin{Definition} {\rm(\cite{gmr.2,r})} Suppose the set $B \subset \left\{X_{t}\right\}_{t \in \mathbb{R}}$, then a function $\psi(\cdot, \cdot)$ defined on $X_{t} \times X_{t}$ is said to be a contractive function on $B \times B$, if  for any sequence $\left\{x_{n}\right\}_{n=1}^{\infty} \subset B$, there is a subsequence $\left\{x_{n_{k}}\right\}_{k=1}^{\infty} \subset\left\{x_{n}\right\}_{n=1}^{\infty}$ such that
$$
\lim _{k \rightarrow \infty} \lim _{l \rightarrow \infty} \psi\left(x_{n_{k}}, x_{n_{l}}\right)=0.
$$
For simplicity, we denote the set of all contractive functions on $B \times B$ by ${\mathcal C}(\widehat{B})$.
\label{def2.6-4}
\end{Definition}

\begin{Lemma} {\rm(\cite{r})}\label{lem2.1-4} Assume the process $\{U(t, \tau)\}_{t \geq \tau}$ has a pullback $\mathcal{D}$-absorbing set $\widehat{B}=\{B(t) : t \in \mathbb{R}\}$ on $\left\{X_{t}\right\}_{t \in \mathbb{R}}$ and there exist $T=T(t, \widehat{B}, \tilde c)=$ $t-\tau$ and $\psi_{t, T}(\cdot, \cdot) \in \operatorname{\mathcal C}(\widehat{B})$ such that
$$
\|U(t, t-T) x-U(t, t-T) y\|_{X_{t}} \leq \tilde c+\psi_{t, T}(x, y),
$$
for any $x, y \in B(\tau)$ and $\tilde c>0$, then $\{U(t, \tau)\}_{t \geq \tau}$ is pullback $\mathcal{D}$-asymptotically compact on $\left\{X_{t}\right\}_{t \in \mathbb{R}}$.

\end{Lemma}

\begin{Lemma} {\rm(\cite{r})}\label{lem2.2-4}
The process $\{U(t, \tau)\}_{t \geq \tau}$ has a time-dependent pullback $\mathcal{D}$-attractor in $\left\{X_{t}\right\}_{t \in \mathbb{R}}$, if it satisfies the following conditions:

(i) $\{U(t, \tau)\}_{t \geq \tau}$ has a pullback $\mathcal{D}$-absorbing set $\widehat{B}_{0}$ in $X_{t}$\\
and

(ii) $\{U(t, \tau)\}_{t \geq \tau}$ is pullback $\mathcal{D}$-asymptotically compact in $\widehat{B}_{0}$.
\end{Lemma}

\section{\large Existence and uniqueness of weak solutions}
To study existence of solutions by the time-dependent pullback $\mathcal D_{\eta_{1}}$-attractors, as in \cite{qy.2}, we first prove the existence and uniqueness theorems of weak solutions to problem $(\ref{1.1-4})$ in this section.

\begin{Definition}
A weak solution of problem $(\ref{1.1-4})$ is a function $u \in C([\tau-k, T]; \mathcal{H}_{t}(\Omega)) \cap L^{2}(\tau, T ; \mathcal{H}_{t}(\Omega)) \cap L^{p}(\tau, T ; L^{p}(\Omega))$ for any $\tau < T \in \mathbb R$, with $u(t)=\phi(t-\tau)$ for all $t \in[\tau-k, \tau]$ such that
\begin{equation}
\begin{aligned}
& \frac{d}{d t}[(u(t), \varphi)+\varepsilon(t)(\nabla u(t), \nabla \varphi)]+(2 a(l(u))-\varepsilon^{\prime}(t))(\nabla u(t), \nabla \varphi) \\
 & = 2(f(u(t)), \varphi)+2(g(t, u_{t}), \varphi)+2(h(t),\varphi),
\end{aligned}
\label{3.1-4}
\end{equation}
for any test function $\varphi \in H_{0}^{1}(\Omega)$.
\label{def3.1-4}
\end{Definition}
\begin{Remark}
The equation $(\ref{3.1-4})$ should be understood in the sense of the generalized function space
$ \mathcal{D}^{\prime}(\tau,+\infty)$.
\end{Remark}
\begin{Corollary}
If $u$ is a weak solution of problem $(\ref{1.1-4})$, then the following energy equality holds
\begin{equation}
\begin{array}{l}
\|u(t)\|^{2}+\varepsilon(t)\|\nabla u(t)\|^{2}+\int_{s}^{t}\left(2 a(l(u))-\varepsilon^{\prime}(r)\right)\|\nabla u(r)\|^{2} d r \\
=\|u(s)\|^{2}+\varepsilon(s)\|\nabla u(s)\|^{2}+2 \int_{s}^{t}(f(u(r))+g(t, u_{r})+h(r), u(r)) dr,
\end{array}
\label{3.2-4}
\end{equation}
for all $\tau  \le s \le t$.
\end{Corollary}
\qquad Firstly, we will use the the standard Faedo-Galerkin method to prove the following theorem.

\begin{Theorem}
Assume that $a(\cdot)$ is a local Lipschitz continuous function and satisfies $(\ref{1.4-4})$, $l(\, \cdot \,)$ is given in $(\ref{1.5-4})$, $f \in C^{1}(\mathbb{R}, \mathbb{R})$ and satisfies $(\ref{1.6-4})-(\ref{1.7-4})$, $h \in L_{l o c}^{2}(\mathbb{R} ; H^{-1}(\Omega))$ and $\phi \in C_{\mathcal H_{t}(\Omega)}$ is given, then for any $\tau \leq t \in \mathbb R$, there exists a weak solution to problem $(\ref{1.1-4})$.
\label{th3.1-4}
\end{Theorem}
$\mathbf{Proof.}$ Suppose the approximate solution $u_{i}(t, x)=\sum\limits_{j=1}^{i} r_{i, j}(t) \omega_{j}(x)$, where $i, j \in {\mathbb N^+}$, $\left\{ {{\omega_j}} \right\}_{j = 1}^\infty $ is a basis of $H^{2}(\Omega) \cap H_{0}^{1}(\Omega)$ and orthonormal in $L^{2}(\Omega)$. The Faedo-Galerkin method needs to find an approximate sequence $\{ {u_i}\}_{i \ge n}$ such that the following approximate system holds
\begin{equation}
\left\{ {\begin{array}{*{20}{l}}
{\frac{d}{{dt}}[({u_i}(t),{\omega _j}) + \varepsilon (t)(\nabla {u_i}(t),\nabla {\omega _j})] + ( {2a(l({u_i})) - {\varepsilon ^\prime }(t)} )(\nabla {u_i}(t),\nabla {\omega _j})}\\
{ = 2(f({u_i}(t)),{\omega _j}) +2(g(t, u_{t,i}),{\omega _j})+ 2\left({h(t),{\omega _j}} \right), \quad \forall\,\, t \in [\tau , + \infty ),}\\
{u_{\tau,i}=\phi,}
\end{array}} \right.
\label{3.3-4}
\end{equation}
where $\phi \in C([-k, 0]; \operatorname{span}\{w_{j}\}_{j=1}^{n}), u_{t, i}=u_{i}(t+\theta)$ and $u_{\tau,i}=u_{i}(\tau+\theta)$ with $\theta \in [-k,0]$.

\textbf{\textbf{\emph{Step}}\,\emph{1:}\,\emph{(A priori estimate for $\boldsymbol{u}$)}}
(1) When $\varepsilon(t)$ is a decreasing function.

Multiplying $(\ref{3.3-4})_{1}$ by the test function ${\gamma _{i,{\rm{ }}j}}(t)$ and then summing $j$ from $1$ to $i$, we obtain
\begin{equation}
\begin{aligned}
&\frac{d}{d t}(\left\|u_{i}(t)\right\|^{2}+\varepsilon(t)\left\|\nabla u_{i}(t)\right\|^{2})+(2a(l(u_{i}))-\varepsilon^{\prime}(t))\left\|\nabla u_{i}(t)\right\|^{2}\\
&= 2\left(f\left(u_{i}(t)\right), u_{i}(t)\right)+2(g(t, u_{t,i}),u_{i}(t))+2\left(h(t), u_{i}(t)\right).
\end{aligned}
\label{3.4-4}
\end{equation}

From (\ref{1.7-4}), we conclude
\begin{equation}
2\left(f\left(u_{i}(t)\right), u_{i}(t)\right) \leq 2 C_{0}|\Omega|-2 C_{2}\left\|u_{i}(t)\right\|_{p}^{p}.
\label{3.5-4}
\end{equation}

Besides, by the Young inequality and assumptions (A1)$-$(A3), we derive
\begin{equation}
2\left(g\left(t, u_{t, i}\right), u_{i}(t)\right) \leq 2 C_{g}\left\|u_{t, i}(t)\right\|_{C_{L^{2}(\Omega)}}^{2}+\|u_{i}(t)\|_{2}^{2} \,.
\label{3.6-4}
\end{equation}

Using the Young and the Cauchy inequalities, it follows that
\begin{equation}
2\left(h(t), u_{i}(t)\right) \leq\|h(t)\|_{-1}^{2}+\left\|\nabla u_{i}(t)\right\|^{2}.
\label{3.7-4}
\end{equation}

Substituting (\ref{3.5-4})$-$(\ref{3.7-4}) into (\ref{3.4-4}), we arrive at
\begin{equation}
\begin{aligned}
&\frac{d}{d t}(\|u_{i}(t)\|^{2}+\varepsilon(t)\|\nabla u_{i}(t)\|^{2})+(2 a(l(u_{i}))-\varepsilon^{\prime}(t))\|\nabla u_{i}(t)\|^{2}+2 C_{2}\|u_{i}(t)\|_{p}^{p} \\
&\leq 2C_{0}|\Omega|+2 C_{g}\|u_{t, i}\|^{2}_{C_{L^{2}(\Omega)}}+\|u_{i}(t)\|^{2}+\|\nabla u_{i}(t)\|^{2}+\|h(t)\|^{2}_{-1}.
\end{aligned}
\label{3.8-4}
\end{equation}

Then from (\ref{1.3-4}) and (\ref{1.4-4}), it follows
\begin{equation}
\begin{aligned}
&\frac{d}{d t}(\|u_{i}(t)\|^{2}+\varepsilon(t)\|\nabla u_{i}(t)\|^{2})+(1+L)\|\nabla u_{i}(t)\|^{2}+2 C_{2}\|u_{i}(t)\|_{p}^{p} \\
&\leq 2C_{0}|\Omega|+2 C_{g}\|u_{t, i}\|^{2}_{C_{L^{2}(\Omega)}}+\|h(t)\|^{2}_{-1}.
\end{aligned}
\label{3.9-4}
\end{equation}

Integrating (\ref{3.9-4}) from $\tau$ to $t$, we deduce
\begin{equation}
\begin{aligned}
&\left\|u_{i}(t)\right\|^{2}+\varepsilon(t)\left\|\nabla u_{i}(t)\right\|^{2}+(1+L)\int_{\tau}^{t}\left\|\nabla u_{i}(s)\right\|^{2} d s+2C_{2}\int_{\tau}^{t}\|u_{i}(s)\|^{p}_{p}ds\\
&\leq 2C_{0}|\Omega|(t-\tau)+\|u_{i}(\tau)\|^{2}+\varepsilon(\tau)\|\nabla u_{i}(\tau)\|^{2}+2C_{g}\int_{\tau}^{t}\|u_{s,i}\|^{2}_{C_{L^{2}(\Omega)}}ds+ \int_{\tau}^{t}\|h(s)\|^{2}_{-1} d s.
\end{aligned}
\label{3.10-4}
\end{equation}

Putting $t+\theta$ instead of $t$ with $\theta \in [-k,0]$ in (\ref{3.10-4}), we obtain
\begin{equation}
\begin{aligned}
&\left\|u_{t, i}\right\|_{C_{L^{2}(\Omega)}}^{2}+\varepsilon_{t}\left\|\nabla u_{t, i}\right\|_{C_{L^{2}(\Omega)}}^{2}+(1+L) \int_{\tau}^{t}\left\|\nabla u_{s, i}\right\|_{C_{L^{2}(\Omega)}}^{2} d s+2 C_{2} \int_{\tau}^{t}\left\|u_{s, i}\right\|_{p}^{p} d s\\
&\leq\|\phi\|_{C_{L^{2}(\Omega)}}^{2}+\varepsilon_{\tau}\|\nabla \phi\|_{C_{L^{2}(\Omega)}}^{2}+2C_{0}|\Omega|(t-\tau)+2 C_{g} \int_{\tau}^{t}\left\|u_{s, i}\right\|_{C_{L^{2}(\Omega)}}^{2} d s+\int_{\tau}^{t}\|h(s)\|^{2}_{-1} d s \text {. }
\end{aligned}
\label{3.11-4}
\end{equation}

By simple calculations, we arrive at
\begin{equation}
\begin{aligned}
&\left\|u_{t, i}\right\|_{C_{L^{2}(\Omega)}^{2}}^{2}+\varepsilon_{t}\left\|\nabla u_{t, i}\right\|_{C_{L^{2}(\Omega)}^{2}}^{2}\leq\|\phi\|_{C_{L^{2}(\Omega)}^{2}}^{2}+\varepsilon_{\tau}\|\nabla \phi\|_{C_{L^{2}(\Omega)}}^{2}+2 C_{0}|\Omega|(t-\tau)\\
&\qquad\qquad\qquad\qquad \qquad+2 C_{g} \int_{\tau}^{t}\left(\left\|u_{s, i}\right\|_{C_{L^{2}(\Omega)}}^{2}+\varepsilon_{s}\left\|\nabla u_{s, i}\right\|_{C_{L^{2}(\Omega)}}^{2}\right) d s+\int_{\tau}^{t}\|h(s)\|_{-1}^{2} d s.
\end{aligned}
\label{3.112-4}
\end{equation}

Thanks to the Gronwall inequality, we conclude
\begin{equation}
\begin{aligned}
&\left\|u_{t, i}\right\|_{C_{L^{2}(\Omega)}}^{2}+\varepsilon_{t}\left\|\nabla u_{t, i}\right\|_{C_{L^{2}(\Omega)}}^{2}\leq e^{2C_{g}(t-\tau)}\left(\|\phi\|_{C_{L^{2}(\Omega)}}^{2}+\varepsilon_{\tau}\|\nabla \phi\|_{C_{L^{2}(\Omega)}}^{2}\right)\\
&\qquad\qquad\qquad\qquad\qquad\qquad\quad+2e^{2C_{g}(t-\tau)}C_{0}|\Omega|(t-\tau)+e^{2C_{g}(t-\tau)}\int_{\tau}^{t}\|h(s)\|^{2}_{-1} d s.
\end{aligned}
\label{3.12-4}
\end{equation}

Then from (\ref{3.11-4}) and (\ref{3.12-4}), we derive
\begin{equation}
\left\{u_{i}\right\} \text { is bounded in } C([\tau-k, T]; \mathcal{H}_{t}(\Omega)) \cap L^{2}(\tau, T ; {H}_{0}^{1}(\Omega)) \cap L^{p}(\tau, T ; L^{p}(\Omega)).
\label{3.13-4}
\end{equation}

Moreover, from (\ref{1.7-4}) we arrive at
\begin{equation}
f(u_{i}(t)) \text { is bounded in } L^{q}(\tau, T ; L^{q}(\Omega)), \text { for all } T>\tau,
\label{3.14-4}
\end{equation}
where $q=\frac{p}{p-1}$ with $p \ge 2$.

Multiplying $(\ref{3.3-4})_{1}$ by $\gamma_{i, j}^{\prime}(t)$ and summing $j$ from 1 to $i$, then by (\ref{1.4-4}), we obtain
\begin{equation}
\begin{aligned}
2\|u_{i}^{\prime}(t)\|^{2}+2\varepsilon(t)\|\nabla u_{i}^{\prime}(t)\|^{2}+m \frac{d}{d t}\left\|\nabla u_{i}(t)\right\|^{2} \leq 2\left(f(u_{i}(t))+g(t, u_{t, i})+h(t), u_{i}^{\prime}(t)\right).
\end{aligned}
\label{3.15-4}
\end{equation}

Furthermore, by $\mathcal F(u)=\int_{0}^{u} f(r) d r$, we arrive at
\begin{equation}
2(f(u_{i}(t)), u_{i}^{\prime}(t))=2 \frac{d}{d t} \int_{\Omega} \mathcal F(u_{i}(t, x)) d x.
\label{3.16-4}
\end{equation}

Using assumptions (A1)$-$(A3) and the Young inequality, we derive
\begin{equation}
2(g(t, u_{t, i}), u_{i}^{\prime}(t)) \leq C_{g}\|u_{t, i}\|_{C_{L^{2}(\Omega)}}^{2}+\|u_{i}^{\prime}(t)\|^{2}.
\label{3.17-4}
\end{equation}

From the Cauchy and Young inequalities, we conclude
\begin{equation}
2(h(t), u_{i}^{\prime}(t)) \leq\|h(t)\|_{-1}^{2}+\|\nabla u_{i}^{\prime}(t)\|^{2} .
\label{3.18-4}
\end{equation}

Substituting (\ref{3.16-4}) $-$ (\ref{3.18-4}) into (\ref{3.15-4}), and by $2\varepsilon(t)-1>\varepsilon(t)$ we obtain
\begin{equation}
\begin{aligned}
&\|u_{i}^{\prime}(t)\|^{2}+\varepsilon(t)\|\nabla u_{i}^{\prime}(t)\|^{2}+m \frac{d}{d t}\|\nabla u_{i}(t)\|^{2} \\
&\leq 2 \frac{d}{d t} \int_{\Omega}\mathcal  F(u_{i}(t, x)) d x+C_{g}\|u_{t, i}\|_{C_{L^{2}(\Omega)}}^{2}+\|h(t)\|_{-1}^{2}.
\end{aligned}
\label{3.19-4}
\end{equation}

Integrating (\ref{3.19-4}) from $\tau$ to $t$, we deduce
\begin{equation}
\begin{aligned}
&\|u_{i}(t)\|^{2}+\int_{\tau}^{t} \varepsilon(s)\left\|u_{i}^{\prime}(s)\right\|^{2} d s+m \int_{\tau}^{t}\|\nabla u_{i}(s)\| \|\nabla u_{i}^{\prime}(s)\| d s+2 \int_{\Omega} \mathcal F(u_{i}(\tau, x)) d x \\
&\leq\|u_{i}(\tau)\|^{2}+2 \int_{\Omega} \mathcal F(u_{i}(t, x)) d x+C_{g} \int_{\tau}^{t}\|u_{s, i}\|_{C_{L^{2}(\Omega)}}^{2} d s+\int_{\tau}^{t}\|h(s)\|^{2}_{-1} d s.
\end{aligned}
\label{3.20-4}
\end{equation}

By (\ref{1.8-4}), it follows that
\begin{equation}
-\widetilde{C}_{0}|\Omega|-\widetilde{C}_{1}\left\|u_{i}(t)\right\|_{p}^{p} \leq 2 \int_{\Omega} \mathcal F\left(u_{i}(t, x)\right) d x \leq 2 \widetilde{C}_{0}-\widetilde{C}_{2}\left\|u_{i}(t)\right\|_{p}^{p}.
\label{3.21-4}
\end{equation}

Substituting (\ref{3.21-4}) into (\ref{3.20-4}), we obtain
\begin{equation}
\begin{aligned}
&\|u_{i}(t)\|^{2}+\int _{\tau}^{t} \varepsilon(s)\|\nabla u_{i}^{\prime}(s)\|^{2} d s+m \int_{\tau}^{t}\left\|\nabla u_{i}^{\prime}(t)\right\|^{2} d s+2 \widetilde{C}_{2}\left\|u_{i}(t)\right\|_{p}^{p}\\
&\leq\|u_{i}(\tau)\|^{2}+4 \widetilde{C}_{0}|\Omega|+2 \widetilde{C}_{1}\|u_{i}(\tau)\|_{p}^{p}+C_{g} \int_{\tau}^{t}\| u_{s,i}\|_{C_{L^{2}(\Omega)}}^{2} d s+\int_{\tau}^{t}\| h(s) \|^{2}_{-1} d s.
\end{aligned}
\label{3.22-4}
\end{equation}

(2) When $\varepsilon(t)$ is an increasing function.

From (\ref{1.3-4}), by some similar calculations to the case (1), we derive
\begin{equation}
\|u_{i}(t)\|^{2}+\int_{\tau}^{t}|\varepsilon(s)|\|\nabla u_{i}^{\prime}(s)\|^{2} d s \leq \widetilde{C}\left(\left\|u_{i}(\tau)\right\|_{p}^{p}+\int_{\tau}^{t}\left\|u_{s, i}\right\|_{C_{L^{2}(\Omega)}}^{2} d s+\int_{\tau}^{t}\|h(s)\|_{-1}^{2} d s\right),
\label{3.23-4}
\end{equation}
where $\widetilde{C} > 0$ is a constant depends on $\widetilde{C}_{0}$, $\widetilde{C}_{1}$ and $C_{g}$.

Putting $t+\theta$ to instead of $t$ with $\theta \in [-k, 0]$ in $(\ref{3.22-4})$ $-$ $(\ref{3.23-4})$ and using the Gronwall inequality, then through similar calculations and estimations to $(\ref{3.11-4})$ and (\ref{3.12-4}), we conclude
\begin{equation}
\left\{u_{i}\right\} \text { is bounded in } L^{\infty}(\tau, T ; H_{0}^{1}(\Omega) \cap L^{p}(\Omega))
\label{3.24-4}
\end{equation}
and
\begin{equation}
\{\partial_{t} u_{i}\} \text { is bounded in } L^{2}(\tau,T ; \mathcal H_{t}(\Omega)).
\label{3.25-4}
\end{equation}

Then from (\ref{3.13-4}), (\ref{3.14-4}), (\ref{3.24-4}), (\ref{3.25-4}), the compactness arguments and the Aubin-Lions lemma (see \cite{Lions}), we derive that there exists a subsequence of $\left\{u_{i}\right\}$ (still marked as $\left\{u_{i}\right\}$), $u \in L^{\infty}\left(\tau, T ; \mathcal{H}_{t}(\Omega)\right) \cap L^{2}(\tau, T ; H_{0}^{1}(\Omega))\cap L^{p}(\tau, T ; L^{p}(\Omega))$ and $\partial_{t} u \in L^{\infty}\left(\tau, T ; \mathcal{H}_{t}(\Omega)\right)$ such that
\begin{equation}
u_{i} \rightharpoonup u \quad \text { weakly-star in } L^{\infty}(\tau-k, T ;\mathcal H_{t}(\Omega));
\label{3.26-4}
\end{equation}
\begin{equation}
u_{i} \rightharpoonup u \quad \text { weakly in } L^{2}(\tau, T ; H_0^1(\Omega));
\label{3.27-4}
\end{equation}
\begin{equation}
u_{i} \rightharpoonup u \quad \text { weakly in } L^{p}(\tau, T ; L^{p}(\Omega));
\label{3.28-4}
\end{equation}
\begin{equation}
f(u_{i}) \rightharpoonup f(u) \quad \text { weakly in } L^{q}(\tau, T ; L^{q}(\Omega));
\label{3.29-4}
\end{equation}
\begin{equation}
a\left(l(u_{i})\right) u_{i} \rightharpoonup a\left(l(u)\right) u \quad \text { weakly in } L^{2}(\tau, T ; H_{0}^{1}(\Omega));
\label{3.30-4}
\end{equation}
\begin{equation}
{\partial _t}{u_i} \rightharpoonup {\partial _t}u \quad \text { weakly in } L^{2}(\tau, T ;\mathcal H_{t}(\Omega));
\label{3.31-4}
\end{equation}
\begin{equation}
u_{i} \rightarrow u \quad \text{ in } C([\tau-k, T];\mathcal H_{t}(\Omega)).
\label{3.32-4}
\end{equation}

From $(\ref{3.26-4})-(\ref{3.32-4})$, we can pass the limit in equation (\ref{3.3-4}) and notice that $\left\{w_{j}\right\}_{j=1}^{\infty}$ is dense in $H_{0}^{1}(\Omega) \cap L^{p}(\Omega)$, (\ref{3.1-4}) holds for any $\varphi \in H_{0}^{1}(\Omega) \cap L^{p}(\Omega)$.

\textbf{\textbf{\emph{Step}}\,\emph{2:}\,\emph{(Verify the initial value)}} In order to verify $u$ is a weak solution to problem (\ref{1.1-4}), we only need to check that $u_{\tau, i}=\phi$.

Choosing a function $\varphi \in C^{1}([\tau, T] ; H_{0}^{1}(\Omega))$ with $\varphi(T)=0$, then we obtain
\begin{equation}
\begin{aligned}
&\int_{\tau}^{T}-(u, \varphi^{\prime}) d s+\int_{\tau}^{T} \int_{\Omega} \varepsilon(s) \nabla(\partial_{s} u) \nabla \varphi d x d s-\int_{\tau}^{T} \int_{\Omega} a(l(u)) (\Delta u)\varphi d x d s \\
&-\int_{\tau}^{T} \int_{\Omega}\left(f(u)+g\left(t, u_{s}\right)+h(s)\right) \varphi d x d s=(u(\tau+\theta), \varphi(\theta)).
\end{aligned}
\label{3.33-4}
\end{equation}

Repeating the same process in the Faedo-Galerkin approximations yields
\begin{equation}
\begin{aligned}
&\int_{\tau}^{T}-(u_{i}, \varphi^{\prime}) d s+\int_{\tau}^{T} \int_{\Omega} \varepsilon(s) \nabla(\partial_{s} u_{i}) \nabla \varphi d x d s-\int_{\tau}^{T} \int_{\Omega} a(l(u_{i})) (\Delta u_{i})\varphi d x d s \\
&-\int_{\tau}^{T} \int_{\Omega}\left(f(u_{i})+g\left(t, u_{s,i}\right)+h(s)\right) \varphi d x d s=(u_{i}(\tau+\theta), \varphi(\theta)).
\end{aligned}
\label{3.34-4}
\end{equation}

Taking limits as $i \to \infty $ in (\ref{3.34-4}) and since $u_{i}(\tau+\theta) \rightarrow \phi(\tau+\theta)$, we deduce
\begin{equation}
\begin{aligned}
&\int_{\tau}^{T}-(u, \varphi^{\prime}) d s+\int_{\tau}^{T} \int_{\Omega} \varepsilon(s) \nabla(\partial_{s} u) \nabla \varphi d x d s-\int_{\tau}^{T} \int_{\Omega} a(l(u)) (\Delta u)\varphi d x d s \\
&-\int_{\tau}^{T} \int_{\Omega}\left(f(u)+g\left(t, u_{s}\right)+h(s)\right) \varphi d x d s=(\phi(\tau+\theta), \varphi(\theta)).
\end{aligned}
\label{3.35-4}
\end{equation}

Then we derive $u(\tau+\theta)=\phi(\theta)$.

Finally, combining the estimations in Step 1 and Step 2, then $u$ is a weak solution of problem $(\ref{1.1-4})$ directly holds. $\hfill$$\Box$

\begin{Theorem}
Under the assumptions of Theorem ${\ref{th3.1-4}}$, if the weak solution of problem $(\ref{1.1-4})$ exists, then it is unique. Moreover, the weak solution depends continuously on its initial value.
\label{th3.2-4}
\end{Theorem}
$\mathbf{Proof.}$ Assuming that $u_{1}$ and $u_{2}$ are two solutions corresponding to the initial values $u_{1}(\tau+\theta)$ and $u_{2}(\tau+\theta)$, respectively, and satisfying
\begin{equation}
\left\{\begin{array}{ll}
\partial_{t}u_{1}-\varepsilon(t) \Delta \partial_{t}u_{1}-a(l(u_{1})) \Delta u_{1}=f(u_{1})+g(t,u_{t,1})+h(t) & \text { in } \Omega \times(\tau, \infty), \\
u_{1}(x,t)=0 & \text { on } \partial \Omega\times(\tau, \infty), \\
u_{1}(x, \tau+\theta)=\phi_{1}(x,\theta),  &\,\, x \in \Omega, \theta\in[-k,0],
\end{array}\right.
\label{3.36-4}
\end{equation}
and
\begin{equation}
\left\{\begin{array}{ll}
\partial_{t}u_{2}-\varepsilon(t) \Delta \partial_{t}u_{2}-a(l(u_{2})) \Delta u_{2}=f(u_{2})+g(t,u_{t,2})+h(t) & \text { in } \Omega \times(\tau, \infty), \\
u_{2}(x,t)=0 & \text { on } \partial \Omega\times(\tau, \infty), \\
u_{2}(x, \tau+\theta)=\phi_{2}(x,\theta),  &\,\, x \in \Omega, \theta\in[-k,0].
\end{array}\right.
\label{3.37-4}
\end{equation}

When $\varepsilon(t)$ is a decreasing function, subtracting $(\ref{3.37-4})_{1}$ from $(\ref{3.36-4})_{1}$, assuming $u=u_{1}-u_{2}$ and taking $L^2$-inner product between $u$ and the resulting equation, we derive
\begin{equation}
\begin{aligned}
&\frac{d}{d t}[\|u\|^{2}+\varepsilon(t)\|\nabla u\|^{2}]+2(a(l(u_{1}))-\varepsilon^{\prime}(t))\|\nabla u\|^{2} \\
&=2(a(l(u_{2}))-a(l(u_{1})))(\nabla u_{2}, \nabla u)+2(f(u_{1})-f(u_{2}), u)+2(g(t,u_{t,1})-g(t,u_{t,2}),u).
\end{aligned}
\label{3.38-4}
\end{equation}

From $(\ref{1.2-4})-(\ref{1.7-4})$, assumptions (A1)$-$(A3) and the Poincar\'{e} inequality, we conclude
\begin{equation}
\frac{d}{d t}(\|u\|^{2}+\varepsilon(t)\|\nabla u \|^{2}) \leq C(\|u_{t}\||_{C_{L^{2}(\Omega)}}^{2}+\|\nabla u_{t}\|_{C_{L^{2}(\Omega)}}^{2}),
\label{3.39-4}
\end{equation}
where $C>0$ is a constant depends on $L$, $m$, $\varepsilon^{\prime}(t)$, $\tilde{\eta}$, $C_{0}, C_{1}$ and $C_{2}$.

Integrating it in $[\tau, t]$ and putting $t+\theta$ instead of $t$, we obtain
\begin{equation}
\|u_{t}\|_{C_{L^{2}(\Omega)}}^{2}+\varepsilon_{t}\|\nabla u_{t}\|_{C_{L^{2}(\Omega)}}^{2} \leq \|u_{\tau}\|_{C_{\mathcal H_{t}(\Omega)}}^{2}+C \int_{\tau}^{t}\left(\left\|u_{s}\right\|_{C_{L^{2}(\Omega)}}^{2}+\|\nabla u_{s}\|_{C_{L^{2}(\Omega)}}^{2}\right) d s.
\label{3.40-4}
\end{equation}

By the Gronwall inequality, we deduce
\begin{equation}
\|u_{t}\|_{C_{L^{2}(\Omega)}}^{2}+\varepsilon_{t}\|\nabla u_{t}\|_{C_{L^{2}(\Omega)}}^{2} \leq e^{C(t+k-\tau)}(\|u_{\tau}\|_{C_{L^{2}(\Omega)}}^{2}+\varepsilon_{\tau}\|\nabla u_{\tau}\|_{C_{L^{2}(\Omega)}}^{2}).
\label{3.41-4}
\end{equation}

When $\varepsilon(t)$ is an increasing function, the following estimates can be obtained from (\ref{1.3-4}), (\ref{1.4-4}) and similar calculations to those above
\begin{equation}
\|u_{t}\|_{C_{L^{2}(\Omega)}}^{2}+|\varepsilon_{t}|\|\nabla u_{t}\|_{C_{L^{2}(\Omega)}}^{2} \leq e^{C(t+k-\tau)}(\|u_{\tau}\|_{C_{L^{2}(\Omega)}}^{2}+\varepsilon_{\tau}\|\nabla u_{\tau}\|_{C_{L^{2}(\Omega)}}^{2}).
\label{3.42-4}
\end{equation}

Consequently, the uniqueness and continuity follow readily.$\hfill$$\Box$

\begin{Corollary}\label{Cor3.1-4}
Thanks to Theorems $\rm{\ref{th3.1-4}}$ and $\rm{\ref{th3.2-4}}$, problem $(\ref{1.1-4})$ has a continuous process
$$
U(t, \tau):C_{\mathcal{ H}_{t}(\Omega)} \rightarrow C_{\mathcal{ H}_{t}(\Omega)}
$$
with $U(t, \tau+\theta)\phi=u(t)$ being the unique weak solution.
\label{re3.1-2}
\end{Corollary}

\section{\large Existence of the time-dependent pullback $\mathcal D_{\eta_{1}}$-attractors}
In this section, we will discuss the existence of the time-dependent pullback $\mathcal D_{\eta_{1}}$-attractor for the process  $\left\{ {U(t,\tau ){\} _{t \ge \tau }}} \right.$ in $C_{\mathcal H_{t}(\Omega)}$ by the similar methods to Zhu and Sun \cite{zs2}. To obtain the pullback $\mathcal D_{\eta_{1}}$-absorbing set, we will first prove the following lemma.
\begin{Lemma}\label{lem4.1-4}
Under the assumptions of Theorems ${\ref{th3.1-4}}$ and ${\ref{th3.2-4}}$, suppose $\phi \in C_{\mathcal H_{t}(\Omega)}$ is given, then the weak solution of problem $(\ref{1.1-4})$ satisfies
\begin{equation}
\begin{aligned}
&\left\|u_{t}\right\|_{C_{L^{2}(\Omega)}}^{2}+\left|\varepsilon_{t}|\|\nabla u_{t}\right\|_{C_{L^{2}(\Omega)}}^{2} \leq \left(1+\frac{C_{g} e^{\eta k}}{\left(1+\lambda_{1}L\right)\eta_{1}}\right) \frac{2 C_{0}|\Omega|}{\eta} e^{\eta k}\\
&\qquad \qquad \qquad \qquad \qquad +\left(1+\frac{C_{g} e^{\eta k}}{\left(1+\lambda_{1}L\right)\left(\eta-\eta_{1}\right)}\right)\left(\|\phi\|_{C_{L^{2}(\Omega)}}^{2}+\varepsilon_{\tau}\|\nabla \phi\|_{C_{L^{2}(\Omega)}}\right) e^{-\eta_{1}(t-\tau)}\\
&\qquad \qquad \qquad \qquad \qquad +\left(1+\frac{C_{g} e^{\eta k}}{\left(1+\lambda_{1}L\right)\left(\eta-\eta_{1}\right)}\right) e^{\eta k} \int_{\tau}^{t} e^{-\eta_{1}(t-s)}\|h(s)\|_{-1}^{2} d s,
\end{aligned}
\label{4.1-4}
\end{equation}
where
$\eta$ and $\eta_{1}$ are constants that satisfy $0<\eta<\frac{(1+L) \lambda_{1}}{1+\lambda_{1} L}$ and $\eta_{1}=\eta-\frac{C_{g}}{1+\lambda_{1} L} e^{\eta k}>0$, respectively.
\end{Lemma}

$\mathbf{Proof.}$ When $\varepsilon(t)$ is a decreasing function, choosing $u$ as the test function of $(\ref{1.1-4})_{1}$ in $L^{2}(\Omega)$, then the weak solution $u$ satisfies
\begin{equation}
\begin{aligned}
&\frac{d}{d t}[\|u(t)\|^{2}+\varepsilon(t)\|\nabla u(t)\|^{2}]+(2 a(l(u))-\varepsilon^{\prime}(t))\|\nabla u\|^{2}=2(f(u)+g(t,u_{t})+h(t), u).
\label{4.2-4}
\end{aligned}
\end{equation}

By (\ref{1.7-4}), we obtain
\begin{equation}
2(f(u), u) \leq 2 C_{0}|\Omega|-2 C_{2}\|u\|_{p}^{p}.
\label{4.3-4}
\end{equation}

Besides, from the Cauchy and Young inequalities, assumptions (A1)$-$(A3), we derive the following inequalities
\begin{equation}
2(g(t, u_{t}), u) \leq C_{g} \|u_{t}\|^{2}_{C_{L^{2}(\Omega)}}+\|u\|^{2}
\label{4.4-4}
\end{equation}
and
\begin{equation}
2(h(t), u) \leq\|h(t)\|_{-1}^{2}+\|\nabla u\|^{2}.
\label{4.5-4}
\end{equation}

Inserting $(\ref{4.3-4})-(\ref{4.5-4})$ into (\ref{4.2-4}) and by the Poincar\'{e} inequality, we deduce
\begin{equation}
\begin{aligned}
&\frac{d}{d t}\left[\|u\|^{2}+\varepsilon(t)\|\nabla u\|^{2}\right]+(2 a(l(u))-\varepsilon^{\prime}(t))\|\nabla u\|^{2}+2C_{2}\|u\|_{p}^{p}\\
&\leq 2C_{0}|\Omega|+C_{g}\|u_{t}\|_{C{_{L^{2}(\Omega)}}} ^{2}+({\lambda^{-1}}+1)\| \nabla u||^{2}+\|h(t)\|_{-1}^{2}.
\end{aligned}
\label{4.6-4}
\end{equation}

Then from (\ref{1.4-4}), we derive
\begin{equation}
\begin{aligned}
&\frac{d}{d t}\left[\|u\|^{2}+\varepsilon(t)\|\nabla u\|^{2}\right]+(1+L)\|\nabla u\|^{2}+2C_{2}\|u\|_{p}^{p}
\leq 2C_{0}|\Omega|+C_{g}\|u_{t}\|_{C{_{L^{2}(\Omega)}}} ^{2}+\|h(t)\|_{-1}^{2}.
\end{aligned}
\label{4.7-4}
\end{equation}

By the Poincar\'{e} inequality, we obtain that there exists $0<\eta<\frac{(1+L) \lambda_{1}}{1+\lambda_{1} L}$ such that
\begin{equation}
\begin{aligned}
&\frac{d}{d t}\left[\|u\|^{2}+\varepsilon(t)\|\nabla u\|^{2}\right]+\eta(\|u\|^{2}+\varepsilon(t)\|\nabla u\|^{2})
\leq 2C_{0}|\Omega|+C_{g}\|u_{t}\|_{C{_{L^{2}(\Omega)}}} ^{2}+\|h(t)\|_{-1}^{2}.
\end{aligned}
\label{4.8-4}
\end{equation}

Multiplying (\ref{4.8-4}) by $e^{\eta t}$ and integrating the resulting inequality from $\tau$ to $t$ yields
\begin{equation}
\begin{aligned}
e^{\eta t}\left(\|u\|^{2}+\varepsilon(t)\|\nabla u\|^{2}\right) & \leq e^{\eta \tau}\left(\|u(\tau)\|^{2}+\varepsilon(\tau)\|\nabla u(\tau)\|^{2}\right)+2C_{0}|\Omega| \int_{\tau}^{t} e^{\eta s} d s\\
&+C_{g} \int_{\tau}^{t} e^{\eta s}\|u_{s}\|_{C_{L^{2}(\Omega)}}^{2} d s+\int_{\tau}^{t} e^{\eta s}\|h(s)\|^{2}_{-1} d s.
\end{aligned}
\label{4.9-4}
\end{equation}

Then putting $t+\theta$ instead of $t$ with $\theta \in [-k,0]$ in (\ref{4.9-4}), we derive
\begin{equation}
\begin{aligned}
&e^{\eta(t+\theta)}\left(\|u(t+\theta)\|^{2}+\varepsilon(t+\theta)\|\nabla u(t+\theta)\|^{2}\right) \leq 2C_{0}|\Omega| \int_{\tau}^{t} e^{\eta s} d s+C_{g} \int_{\tau}^{t} e^{\eta s}\|u_{s}\|_{C_{L^{2}(\Omega)}}^{2} d s\\
&\qquad \qquad \qquad \qquad +e^{\eta(\tau+\theta)}\left(\|u(\tau+\theta)\|^{2}+\varepsilon(\tau+\theta)\|\nabla u(\tau+\theta)\|^{2}\right)
+\int_{\tau}^{t} e^{\eta s}\|h(s)\|_{- 1}^{2} d s.
\end{aligned}
\label{4.10-4}
\end{equation}

Noting that $\theta \in [-k,0]$ with $k>0$, then from (\ref{4.10-4}), it follows
\begin{equation}
\begin{aligned}
e^{\eta t}\left(\|u_{t}\|_{C_{L^{2}(\Omega)}}^{2}+\varepsilon_{t}\|\nabla u_{t}\|^{2}_{C_{L^{2}(\Omega)}}\right) & \leq e^{\eta \tau}\left(\|\phi\|_{C_{L^{2}(\Omega)}}^{2}+\varepsilon_{\tau}\|\nabla \phi\|_{C_{L^{2}(\Omega)}}^{2}\right)+2C_{0}|\Omega| e^{\eta k}\int_{\tau}^{t} e^{\eta s} d s\\
&+C_{g} e^{\eta k} \int_{\tau}^{t} e^{\eta s}\|u_{s}\|_{C_{L^{2}(\Omega)}}^{2} d s+ e^{\eta k}\int_{\tau}^{t} e^{\eta s}\|h(s)\|^{2}_{-1} d s.
\end{aligned}
\label{4.11-4}
\end{equation}

From (\ref{1.3-4}) and the Poincar\'{e} inequality, we obtain
\begin{equation}
\begin{aligned}
e^{\eta t}\|u_{t}\|_{C_{L^{2}(\Omega)}}^{2} & \leq \frac{1}{1+\lambda_{1}L} e^{\eta \tau}\left(\|\phi\|_{C_{L^{2}(\Omega)}}^{2}+\varepsilon_{\tau}\|\nabla \phi\|_{C_{L^{2}(\Omega)}}^{2}\right)+\frac{2C_{0}|\Omega|}{1+\lambda_{1}L}  e^{\eta k}\int_{\tau}^{t} e^{\eta s} d s\\
&+\frac{C_{g}}{1+\lambda_{1}L} e^{\eta k} \int_{\tau}^{t} e^{\eta s}\|u_{s}\|_{C_{L^{2}(\Omega)}}^{2} d s+ \frac{1}{1+\lambda_{1}L}e^{\eta k}\int_{\tau}^{t} e^{\eta s}\|h(s)\|^{2}_{-1} d s.
\end{aligned}
\label{4.12-4}
\end{equation}

Now, we apply the Gronwall lemma to (\ref{4.11-4}).
Let
\begin{equation}
w(t)=e^{\eta t}\left\|u_{t}\right\|_{C_{L^{2}(\Omega)}}^{2}
\label{4.13-4}
\end{equation}
and
\begin{equation}
\begin{aligned}
v(t) & = \frac{1}{1+\lambda_{1}L} e^{\eta \tau}\left(\|\phi\|_{C_{L^{2}(\Omega)}}^{2}+\varepsilon_{\tau}\|\nabla \phi\|_{C_{L^{2}(\Omega)}}^{2}\right)+\frac{2C_{0}|\Omega|}{1+\lambda_{1}L}  e^{\eta k}\int_{\tau}^{t} e^{\eta s} d s\\
&+ \frac{1}{1+\lambda_{1}L}e^{\eta k}\int_{\tau}^{t} e^{\eta s}\|h(s)\|^{2}_{-1} d s.
\end{aligned}
\label{4.14-4}
\end{equation}

Taking $a=\tau$, then by (\ref{4.14-4}), we obtain
\begin{equation}
v(a)=v(\tau)=\frac{1}{1+\lambda_{1}L} e^{\eta \tau}\left(\|\phi\|_{C_{L^{2}(\Omega)}}^{2}+\varepsilon_{\tau}\|\nabla \phi\|_{C_{L^{2}(\Omega)}}^{2}\right).
\label{4.15-4}
\end{equation}

Let $\bar{u}(s)=\frac{C_{g}}{1+\lambda_{1} L} e^{\eta k}$, then we derive
\begin{equation}
e^{\int_{a}^{t} \bar{u}(s) d s}=e^{\frac{C_{g}}{1+\lambda_{1} L} e^{\eta k}\cdot (t-\tau)}
\label{4.16-4}
\end{equation}

and
\begin{equation}
e^{\int_{s}^{t} \bar{u}(r) d r}=e^{\frac{C_{g}}{1+\lambda_{1} L} e^{\eta k}\cdot (t-s)}.
\label{4.17-4}
\end{equation}

Then from (\ref{4.15-4}) and (\ref{4.16-4}), we arrive at
\begin{equation}
v(a) e^{\int_{a}^{t}\bar u(s) d s}=\frac{1}{1+\lambda_{1} L} e^{\eta \tau}\left(\|\phi\|_{C_{L^{2}(\Omega)}}^{2}+\varepsilon_{\tau}\|\nabla \phi\|_{C_{L^{2}(\Omega)}}^{2}\right) e^{\frac{C_{g}}{1+\lambda_{1} L} e^{\eta k} \cdot(t-\tau)}.
\label{4.18-4}
\end{equation}

Besides, by (\ref{4.14-4}), we deduce
\begin{equation}
\frac{d v(s)}{d s}=\frac{2 C_{0}|\Omega|}{1+\lambda_{1} L} e^{\eta k} e^{\eta s}+\frac{1}{1+\lambda_{1} L} e^{\eta k} e^{\eta s}\|h(t)\|_{-1}^{2}.
\label{4.19-4}
\end{equation}

From (\ref{4.17-4}) and (\ref{4.19-4}), we get
\begin{equation}
\begin{aligned}
&\int_{a}^{t}e^{\int_{s}^{t} \bar{u}(r)d r} \cdot \frac{d v}{d s} d s=e^{\frac{C_{g}}{1+\lambda_{1} L} e^{\eta k} \cdot t}\cdot\frac{2 C_{0}|\Omega|}{1+\lambda_{1} L} e^{\eta k} \int_{\tau}^{t} e^{s\left(\eta-\frac{C_{g}}{1+\lambda_{1} L} e^{\eta k}\right)} d s\\
&\qquad\qquad\qquad\qquad\,\,+e^{\frac{C_{g}}{1+\lambda_{1} L} e^{\eta k} \cdot t}\cdot\frac{1}{1+\lambda_{1} L} e^{\eta k}\int_{\tau}^{t} e^{s\left(\eta-\frac{C_{g}}{1+\lambda_{1} L} e^{\eta k}\right)}\|h(s)\|_{-1}^{2} d s.
\label{4.20-4}
\end{aligned}
\end{equation}

Let $\eta_{1}=\eta-\frac{C_{g}}{1+\lambda_{1}L} e^{\eta k}>0$, then by (\ref{4.17-4}), it follows
\begin{equation}
\int_{a}^{t}e^{\int_{s}^{t} \bar{u}(r)d r} \cdot \frac{d v}{d s} d s \leq
\frac{2 C_{0}|\Omega|}{\left(1+\lambda_{1} L\right) \eta_{1}} e^{\eta(k+t)}+\frac{1}{1+\lambda_{1} L} e^{\eta k} e^{\frac{C_{g}}{1+\lambda_{1} L} e^{\eta k} \cdot t} \int_{t}^{t} e^{\eta_{1} s}\|h(s)\|_{-1}^{2} d s.
\label{4.21-4}
\end{equation}

Then by (\ref{4.18-4}), (\ref{4.21-4}) and the Gronwall inequality, we obtain
\begin{equation}
\begin{aligned}
&w(t)=e^{\eta t}\|u_{t}\|_{C_{L^{2}(\Omega)}}^{2} \leq v(a) e^{\int_{a}^{t}\bar u(s) d s}+\int_{a}^{t} e^{\int_{s}^{t} \bar    u(r) d r} \cdot \frac{d v}{d s} d s\\
&\qquad=\frac{1}{1+\lambda_{1} L} e^{\eta \tau}\left(\|\phi\|_{C_{L^{2}(\Omega)}}^{2}+\varepsilon_{\tau}\|\nabla \phi\|_{C_{L^{2}(\Omega)}}^{2}\right) e^{\frac{C_{g}}{1+\lambda_{1} L} e^{\eta k} \cdot(t-\tau)}\\
&\qquad+\frac{2 C_{0}|\Omega|}{\left(1+\lambda_{1} L) \eta_{1}\right.} e^{\eta(k+t)}+\frac{1}{1+\lambda_{1} L} e^{\eta k} e^{\frac{C_{g}}{1+\lambda_{1} L} e^{\eta k} \cdot t} \int_{\tau}^{t} e^{\eta_{1} s}\|h(s)\|^{2}_{-1} d s.
\end{aligned}
\label{4.22-4}
\end{equation}

Dividing both sides of (\ref{4.11-4}) by $e^{\eta t}$, we conclude
\begin{equation}
\begin{aligned}
\|u_{t}\|_{C_{L^{2}(\Omega)}}^{2}+ & \varepsilon_{t}\|\nabla u_{t}\|_{C_{L^{2}(\Omega)}}^{2}\leq e^{\eta (\tau-t)}\left(\|\phi\|_{C_{L^{2}(\Omega)}}^{2}+\varepsilon_{\tau}\|\nabla \phi\|_{C_{L^{2}(\Omega)}}^{2}\right)+2 C_{0}|\Omega| e^{\eta (k-t)} \int_{\tau}^{t} e^{\eta s} d s \\
&\qquad\qquad\qquad+C_{g} e^{\eta (k-t)} \int_{\tau}^{t} e^{\eta s}\left\|u_{s}\right\|_{C_{L^{2}(\Omega)}}^{2} d s+e^{\eta (k-t)} \int_{\tau}^{t} e^{\eta s}\|h(s)\|_{-1}^{2} d s.
\end{aligned}
\label{4.23-4}
\end{equation}

After some simple calculations, we derive the following inequalities
\begin{equation}
(a)^{1}=e^{\eta(\tau-t)}\left(\|\phi\|_{C_{L^{2}(\Omega)}}^{2}+\varepsilon_{\tau}\|\nabla \phi\|_{C_{L^{2}(\Omega)}}^{2}\right) \leq e^{-\eta_{1}(t-\tau)}\left(\|\phi\|_{C_{L^{2}(\Omega)}}^{2}+\varepsilon_{\tau}\|\nabla \phi\|_{C_{L^{2}(\Omega)}}^{2}\right),
\label{4.24-4}
\end{equation}
\begin{equation}
(b)^{1}=2 C_{0}|\Omega| e^{\eta (k-t)} \int_{\tau}^{t} e^{\eta s} d s \leq \frac{2 C_{0} |\Omega|}{\eta} e^{\eta k}
\label{4.25-4}
\end{equation}
and
\begin{equation}
(c)^{1}=e^{\eta (k-t)} \int_{\tau}^{t} e^{\eta s}\|h(s)\|_{-1}^{2} d s \leq e^{\eta k} \int_{\tau}^{t} e^{-\eta_{1}(t-s)}\|h(s)\|_{-1}^{2} d s.
\label{4.26-4}
\end{equation}

Substituting (\ref{4.22-4}), $(\ref{4.24-4})-(\ref{4.26-4})$ into (\ref{4.23-4}), we deduce
\begin{equation}
\begin{aligned}
\|u_{t}\|_{C_{L^{2}(\Omega)}}^{2}+ & \varepsilon_{t}\|\nabla u_{t}\|_{C_{L^{2}(\Omega)}}^{2}\leq (a)^{1}+(b)^{1}+(c)^{1}+(a)+(b)+(c),
\end{aligned}
\label{4.27-4}
\end{equation}
where
\begin{equation}
(a)=C_{g} e^{\eta(k-t)} \int_{\tau}^{t} \frac{1}{1+\lambda_{1} L} e^{\eta \tau}\left(\|\phi\|_{C_{L^{2}(\Omega)}}^{2}+\varepsilon_{\tau}\|\nabla \phi\|_{C_{L^{2}(\Omega)}}^{2}\right) e^{\frac{C_{g}}{1+\lambda_{1} L} e^{\eta k} \cdot(s-\tau)} d s,
\label{4.28-4}
\end{equation}
\begin{equation}
(b)=C_{g} e^{\eta(k-t)} \int_{\tau}^{t} \frac{2 C_{0}\left|\Omega\right|}{\left(1+\lambda_{1} L\right)\eta_{1}} e^{\eta(k+s)} d s
\label{4.29-4}
\end{equation}
and
\begin{equation}
(c)=C_{g} e^{\eta(k-t)} \int_{\tau}^{t} \frac{1}{1+\lambda_{1} L} e^{\eta k} e^{\frac{C_{g}}{1+\lambda_{1} L} e^{\eta k} \cdot s}\left(\int_{\tau}^{s} e^{\eta_{1} s}\|h(s)\|_{-1}^{2} d s\right) d s.
\label{4.30-4}
\end{equation}

Next, we will estimate $(a)$ to $(c)$ in turn. By some simple estimates, we obtain
$$
(a)\leq \frac{C_{g} e^{\eta k}}{\left.\left(1+\lambda_{1} L\right)(\eta-\eta_{1}\right)}\left(\|\phi\|_{C_{L^{2}(\Omega)}}^{2}+\varepsilon_{\tau}\|\nabla \phi\|_{C_{L^{2}(\Omega)}}^{2}\right) e^{-\eta_{1}(t-\tau)},
$$
$$
(b) \leq \frac{C_{g} e^{\eta k}}{\left(1+\lambda_{1} L\right) \eta_{1}}\cdot \frac{2 C_{0}|\Omega|}{\eta} e^{\eta k}
$$
and
$$
(c)\leq \frac{C_{g} e^{\eta k}}{\left(1+\lambda_{1} L\right)\left(\eta-\eta_{1}\right)} e^{\eta k} \int_{\tau}^{t} e^{-\eta_{1}(t-s)} \|h(s)\|^{2}_{-1}d s .
$$

From the above estimates and $(\ref{4.28-4})-(\ref{4.30-4})$, it follows
\begin{equation}
(a)^{1}+(a) \leq\left(1+\frac{C_{g} e^{\eta k}}{\left(1+\lambda_{1}L)\left(\eta-\eta_{1}\right)\right.}\right)\left(\|\phi\|_{C_{L^{2}(\Omega)}}^{2}+\varepsilon_{\tau}\|\nabla \phi\|_{C_{L^{2}(\Omega)}}^{2}\right) e^{-\eta_{1}(t-\tau)},
\label{4.31-4}
\end{equation}
\begin{equation}
(b)^{1}+(b) \leq\left(1+\frac{\left. C_{g} e^{\eta k}\right.}{\left(1+\lambda_{1}L\right)\eta_{1}}\right) \frac{2 C_{0}|\Omega|}{\eta} e^{\eta k}
\label{4.32-4}
\end{equation}
and
\begin{equation}
(c)^{1}+(c) \leq\left(1+\frac{\left. C_{g} e^{\eta k}\right.}{\left(1+\lambda_{1} L)\left(\eta-\eta_{1}\right)\right.}\right) e^{\eta k} \int_{\tau}^{t} e^{-\eta_{1}(t-s)}\|h(s)\|^{2}_{-1}d s.
\label{4.33-4}
\end{equation}

Inserting $(\ref{4.31-4})-(\ref{4.33-4})$ into (\ref{4.27-4}), we conclude
\begin{equation}
\begin{aligned}
&\left\|u_{t}\right\|_{C_{L^{2}(\Omega)}}^{2}+\varepsilon_{t} \|\nabla u_{t}\|_{C_{L^{2}(\Omega)}}^{2} \leq \left(1+\frac{C_{g} e^{\eta k}}{\left(1+\lambda_{1}L\right)\eta_{1}}\right) \frac{2 C_{0}|\Omega|}{\eta} e^{\eta k}\\
&\qquad \qquad \qquad \qquad \qquad +\left(1+\frac{C_{g} e^{\eta k}}{\left(1+\lambda_{1}L\right)\left(\eta-\eta_{1}\right)}\right)\left(\|\phi\|_{C_{L^{2}(\Omega)}}^{2}+\varepsilon_{\tau}\|\nabla \phi\|_{C_{L^{2}(\Omega)}}\right) e^{-\eta_{1}(t-\tau)}\\
&\qquad \qquad \qquad \qquad \qquad +\left(1+\frac{C_{g} e^{\eta k}}{\left(1+\lambda_{1}L\right)\left(\eta-\eta_{1}\right)}\right) e^{\eta k} \int_{\tau}^{t} e^{-\eta_{1}(t-s)}\|h(s)\|_{-1}^{2} d s.
\end{aligned}
\label{4.34-4}
\end{equation}

Finally, we will consider the case when $\varepsilon(t)$ is an increasing function. In fact, (\ref{4.1-4}) can be obtained from (\ref{1.3-4}) and (\ref{1.4-4}) through calculations similar to (\ref{4.34-4}), which we will omit here.$\hfill$$\Box$

In order to obtain the existence of the time-dependent pullback $\mathcal D_{\eta_{1}}$-absorbing set, we further introduce the following concept.

\begin{Definition} (Tempered universe)\label{def4.1-4}
 For each $\eta_{1}>0$, let $\mathcal{D}_{\eta_{1}}$ be the class of all families of nonempty subsets $\widehat{D}=\{D(t): t \in \mathbb{R}\} \subset \Gamma\left(C_{\mathcal H_{t}(\Omega)}\right)$ such that
$$
\lim _{\tau \rightarrow-\infty}\left(e^{\eta_{1} \tau} \sup _{u \in D(\tau)}{\|u\|}^{2}_{C_{\mathcal H_{t}(\Omega)}}\right)=0.
$$
\end{Definition}

\begin{Lemma}\label{lem4.3-4}
Under the assumption of Theorem $\ref{lem4.1-4}$, if $h(t)$ also satisfies
\begin{equation}
\int_{-\infty}^{t} e^{\eta_{1} s}\|h(s)\|^{2}_{-1} d s<\infty
\label{4.35-4}
\end{equation}
for any $ t \in \mathbb{R}$, then the family $\widehat{D}_{1}=\left\{D_{1}(t): t \in \mathbb{R}\right\}$ with $D_{1}(t)=\mathcal{\bar{B}}_{C_{\mathcal H_{t}(\Omega)}}\left(0, {\rho}(t)\right)$, the closed ball in $C_{\mathcal H_{t}(\Omega)}$ of centre zero and radius $\rho(t)$, where
\begin{equation}
\begin{aligned}
\rho^{2}(t) &=C_{3}\left(1+\frac{\left.C_{g} e^{\eta k}\right.}{\left(1+\lambda_{1} L\right) \eta_{1}}\right) \frac{2 C_{0}|\Omega|}{\eta} e^{\eta k} \\
&\left.+\left(1+\frac{C_{g} e^{\eta k}}{\left(1+\lambda_{1} L\right)\left(\eta-\eta_{1}\right)}\right) e^{\eta k} \int_{\tau}^{t} e^{-\eta_{1}(t-s)} \| h(s) \|_{-1}^{2} \right.d s
\end{aligned}
\label{4.36-2}
\end{equation}
and $C_{3}=C(\eta, \eta_{1}, \lambda_{1}, L, k)$, is the time-dependent pullback $\mathcal{D}_{\eta_{1}}$-absorbing family of process ${\{ U(t,\tau )\} _{t \ge \tau }}$ in $C_{\mathcal H_{t}(\Omega)}$. Moreover, $\widehat{D}_{1} \in \mathcal{D}_{\eta_{1}}$.
\label{lem4.2-4}
\end{Lemma}
$\mathbf{Proof.}$  From Lemma \ref{lem4.1-4} and Definition \ref{def4.1-4}, we can derive that $\widehat{D}_{1}$ is the time-dependent pullback $\mathcal D_{\eta_{1}}$-absorbing set. Besides, by (\ref{4.36-2}), we obtain $e^{\eta_{1} t} \rho^{2}(t) \rightarrow 0$ as $t \rightarrow-\infty$. Then thanks to Definition \ref{def4.1-4}, $\widehat{D}_{1} \in \mathcal{D}_{\eta_{1}}$ follows directly. $\hfill$$\Box$

Next, we will use the contractive function method to verify the existence of the time-dependent pullback $\mathcal D_{\eta_{1}}$-attractor for the process $\left\{ {U(t,\tau ){\} _{t \ge \tau }}} \right.$ of problem (\ref{1.1-4}). The following lemma can be proved by similar calculations to Theorems \ref{th3.1-4} and \ref{th3.2-4}.

\begin{Lemma}
Under the assumptions of Theorem $\ref{lem4.2-4}$, assume $\phi \in C_{\mathcal H_{t}(\Omega)}$ is given, if $\left\{u^{s}(t)\right\}_{s \in \mathbb{N}^{+}}$ is a sequence of the solutions to problem $(\ref{1.1-4})$ with initial data $u^{s}(\tau) \in C_{\mathcal H_{t}(\Omega)}$, then there exists a subsequence of $\left\{u^{s}(t)\right\}_{s \in \mathbb{N}^{+}}$ that convergence strongly in $L^{2}(\tau, T;L^{2}(\Omega))$.
\label{lem4.4-4}
\end{Lemma}

Now we will establish the pullback $\mathcal D_{\eta_{1}}$-asymptotic compactness for the process $\left\{ {U(t,\tau ){\} _{t \ge \tau }}} \right.$ of problem (\ref{1.1-4}).

\begin{Lemma}
Under the assumptions of Lemma $\ref{lem4.3-4}$, if $a(\cdot)$ is locally Lipschitz continuous, then the process $\left\{ {U(t,\tau ){\} _{t \ge \tau }}} \right.$ is pullback $\mathcal D_{\eta_{1}}$-asymptotically compact in $C_{\mathcal H_{t}(\Omega)}$.
\label{lem4.3-2}
\end{Lemma}
$\mathbf{Proof.}$ Suppose $u^{j}(t)$ is a weak solution of problem (\ref{1.1-4}) corresponding to the initial value $\phi^{j}(x, \theta) \in D_{1}(\tau)$ $(j=1,2)$, and assume $\tilde{u}=u^{1}-u^{2}$, then we arrive at
\begin{equation}
\partial_{t} \tilde{u}-\varepsilon(t) \partial_{t} \Delta \tilde{u}-a(l(u)) \Delta u^{1}+a(l(u^{2})) \Delta u^{2}=f(u^{1})-f(u^{2})+g(t, u_{t}^{1})-g(t, u_{t}^{2})
\label{4.37-4}
\end{equation}
with initial data
\begin{equation}
\tilde u^{j}(x, \tau+\theta)=\phi^{j}(x, \theta), \quad x \in \Omega, \,\,\theta \in[-k, 0].
\label{4.38-4}
\end{equation}

We shall first discuss the case when $\varepsilon(t)$ is a decreasing function.

Choosing $\tilde u$ as the test function of $(\ref{4.37-4})$, then we obtain
\begin{equation}
\begin{aligned}
&\frac{d}{a t}[\|\tilde{u}\|^{2}+\varepsilon(t)\|\nabla \tilde{u}\|^{2}]+2(a(l(u^{1}))-a(l(u^{2})))\|\nabla \tilde{u}\|^{2} \\
=&2(a(l(u^{2}))-a(l(u^{1})))(\nabla u^{2}, \nabla \tilde{u})+2(f(u^{1})-f(u^{2}), \tilde{u})+2(g(t, u_{t}^{1})-g(t, u_{t}^{2}), \tilde{u}).
\end{aligned}
\label{4.39-4}
\end{equation}

Thanks to (\ref{1.6-4}), we conclude
\begin{equation}
2(f(u^{1})-f(u^{2}), \tilde{u}) \leq 2 \tilde{\eta} \|\tilde{u}\|^{2}.
\label{4.40-4}
\end{equation}

From assumptions (A1)$-$(A3), we deduce
\begin{equation}
2\left(g(t, u_{t}^{1})-g(t, u_{t}^{2}), \tilde{u}\right) \leq 2 C_{g}\|u_{t}^{1}-u_{t}^{2}\|_{C_{L^{2}(\Omega)}}^{2}\|\tilde{u}\|.
\label{4.41-4}
\end{equation}

Noting that $a(\cdot)$ is locally Lipschitz continuous, we can obtain the following inequalities from the Young and Cauchy inequalities
\begin{equation}
2(a(l(u^{2}))-a(l(u^{1})))(\nabla u^{2}, \nabla \tilde{u})\leq 2 m\|\nabla u^{1}-\nabla u^{2}\|^{2}+\frac{(L_{a}(R))^{2}\|l\|^{2}\|\nabla u^{1}\|^{2}\|u^{1}-u^{2}\|^{2}}{2 m},
\label{4.42-4}
\end{equation}
where $L_{a}(R)$ is the Lipschitz constant of the nonlocal term $a(\cdot)$ in $[-R, R]$.

Substituting $(\ref{4.40-4})-(\ref{4.41-4})$ into (\ref{4.39-4}), after some simple calculations, we conclude
\begin{equation}
\frac{d}{d t}[\|\tilde{u}\|^{2}+\varepsilon (t)\|\nabla \tilde{u}\|^{2}]+2 \eta\left(\|\tilde{u}\|^{2}+\varepsilon (t)\|\nabla \tilde{u}\|^{2}\right) \leq 2 \tilde{\eta}\|\tilde{u}\|^{2}+2\left.C_{g}\left\|u_{t}^{1}-u_{t}^{2}\right\|_{C_{L^{2}(\Omega)}}\|\tilde{u}\|\right.,
\label{4.43-4}
\end{equation}
where parameter $\eta>0$ is the same as in Lemma \ref{lem4.1-4}.

Then by the Gronwall inequality, we obtain
\begin{equation}
\begin{aligned}
\|\tilde{u}(t)\|^{2}+\varepsilon(t)\|\nabla \tilde{u}(t)\|^{2} &\leq(\| \tilde{u}(\tau)\|^{2}+\varepsilon(\tau)\| \nabla \tilde{u}(\tau) \|^{2}) e^{- 2\eta(t-\tau)}+2 \tilde{\eta} \int_{\tau}^{t}\|\tilde{u}(s)\|^{2} d s \\
&+2 C_{g} e^{-2 \eta t} \int_{\tau}^{t} e^{2 \eta s}\left\|u_{s}^{1}-u_{s}^{2}\right\|_{C_{L^{2}(\Omega)}}\|\tilde{u}(s)\| d s.
\end{aligned}
\label{4.44-4}
\end{equation}

Putting $t+\theta$ instead of $t$ with $\theta \in [-k,0]$ in (\ref{4.44-4}), we conclude
\begin{equation}
\begin{aligned}
\|\tilde{u}_{t}\|_{C_{L^{2}(\Omega)}}^{2}+&\varepsilon_{t}\|\nabla \tilde{u}(t)\|^{2}_{C_{L^{2}(\Omega)}} \leq\left(\| \tilde{u}_{\tau}\|_{C_{L^{2}(\Omega)}}^{2}+\varepsilon_{\tau}\| \nabla \tilde{u}_{\tau} \|_{C_{L^{2}(\Omega)}}^{2}\right) e^{- 2\eta(t-k-\tau)}\\
&\qquad\qquad \,\, +2 \tilde{\eta} \int_{\tau}^{t}\|\tilde{u}(s)\|^{2} d s +2 C_{g} e^{-2 \eta (t-k)} \int_{\tau}^{t} e^{2 \eta s}\|u_{s}^{1}-u_{s}^{2}\|_{C_{L^{2}(\Omega)}}\|\tilde{u}(s)\| d s.
\end{aligned}
\label{4.45-4}
\end{equation}

Thanks to the H$\rm \ddot{o}$lder inequality, we arrive at
\begin{equation}
\begin{aligned}
&2 C_{g} e^{-2 \eta (t-k)} \int_{\tau}^{t} e^{2 \eta s}\|u_{s}^{1}-u_{s}^{2}\|_{C_{L^{2}(\Omega)}}\|\tilde{u}(s)\| d s\\
&\leq 2 C_{g} e^{-2 \eta (t-k)}\left(\int_{\tau}^{t} e^{4 \eta s}\|u_{s}^{1}-u_{s}^{2}\|_{C_{L^{2}(\Omega)}}^{2} d s\right)^{\frac{1}{2}}\left(\int_{\tau}^{t}\|\tilde{u}(s)\|^{2} d s\right)^{\frac{1}{2}}\\
&\leq 4 C_{g} e^{-2 \eta (t-k)}\left(\int_{\tau}^{t} e^{4 \eta s}(\|u_{s}^{1}\|_{C_{L^{2}(\Omega)}}^{2} +\|u_{s}^{2}\|_{C_{L^{2}(\Omega)}}^{2}) d s\right)^{\frac{1}{2}}\left(\int_{\tau}^{t}\|\tilde{u}(s)\|^{2} d s\right)^{\frac{1}{2}}.
\end{aligned}
\label{4.46-4}
\end{equation}

From $\eta_{1}<\eta$ and (\ref{1.3-4}), we derive
\begin{equation}
\left(\| \tilde{u}_{\tau}\|_{C_{L^{2}(\Omega)}}^{2}+\varepsilon_{\tau}\| \nabla \tilde{u}_{\tau} \|_{C_{L^{2}(\Omega)}}^{2}\right) e^{- 2\eta(t-k-\tau)}\leq\left(\left\|\tilde{u}_{\tau}\right\|_{C_{L^{2}(\Omega)}}^{2}+L\left\|\nabla \tilde{u}_{\tau}\right\|_{C_{L^{2}(\Omega)}}^{2}\right) e^{-2 \eta_{1}(t-k-\tau)}.
\label{4.47-4}
\end{equation}

Substituting (\ref{4.46-4}) and (\ref{4.47-4}) into (\ref{4.45-4}), we deduce
\begin{equation}
\begin{aligned}
&\|\tilde{u}_{t}\|_{C_{L^{2}(\Omega)}}^{2}+\varepsilon_{t}\|\nabla \tilde{u}(t)\|^{2}_{C_{L^{2}(\Omega)}} \\ &\leq\left(\left\|\tilde{u}_{\tau}\right\|_{C_{L^{2}(\Omega)}}^{2}+L\left\|\nabla \tilde{u}_{\tau}\right\|_{C_{L^{2}(\Omega)}}^{2}\right) e^{-2 \eta_{1}(t-k-\tau)}+2 \tilde{\eta} \int_{\tau}^{t}\|\tilde{u}(s)\|^{2} d s\\
&+ 4 C_{g} e^{-2 \eta (t-k)}\left(\int_{\tau}^{t} e^{4 \eta s}(\|u_{s}^{1}\|_{C_{L^{2}(\Omega)}}^{2} +\|u_{s}^{2}\|_{C_{L^{2}(\Omega)}}^{2}) d s\right)^{\frac{1}{2}}\left(\int_{\tau}^{t}\|\tilde{u}(s)\|^{2} d s\right)^{\frac{1}{2}}.
\end{aligned}
\label{4.48-4}
\end{equation}

When $\varepsilon(t)$ is an increasing function, the following inequality can be derived from some estimates like (\ref{4.48-4})
\begin{equation}
\begin{aligned}
&\|\tilde{u}_{t}\|_{C_{L^{2}(\Omega)}}^{2}+|\varepsilon_{t}|\|\nabla \tilde{u}(t)\|^{2}_{C_{L^{2}(\Omega)}} \\ &\leq\left(\left\|\tilde{u}_{\tau}\right\|_{C_{L^{2}(\Omega)}}^{2}+L\left\|\nabla \tilde{u}_{\tau}\right\|_{C_{L^{2}(\Omega)}}^{2}\right) e^{-2 \eta_{1}(t-k-\tau)}+2 \tilde{\eta} \int_{\tau}^{t}\|\tilde{u}(s)\|^{2} d s\\
&+ 4 C_{g} e^{-2 \eta (t-k)}\left(\int_{\tau}^{t} e^{4 \eta s}(\|u_{s}^{1}\|_{C_{L^{2}(\Omega)}}^{2} +\|u_{s}^{2}\|_{C_{L^{2}(\Omega)}}^{2}) d s\right)^{\frac{1}{2}}\left(\int_{\tau}^{t}\|\tilde{u}(s)\|^{2} d s\right)^{\frac{1}{2}}.
\end{aligned}
\label{4.49-4}
\end{equation}

Furthermore, let $T=t-\tau$ and
\begin{equation}
\begin{aligned}
&\psi_{t, T}\left(u^{1}, u^{2}\right)=2 \tilde{\eta} \int_{\tau}^{t}\|\tilde{u}(s)\|^{2} d s+ 4 C_{g} e^{-2 \eta (t-k)}\left(\int_{\tau}^{t} e^{4 \eta s}(\|u_{s}^{1}\|_{C_{L^{2}(\Omega)}}^{2} +\|u_{s}^{2}\|_{C_{L^{2}(\Omega)}}^{2}) d s\right)^{\frac{1}{2}}\\
&\qquad\qquad\quad\,\,\times\left(\int_{\tau}^{t}\|\tilde{u}(s)\|^{2} d s\right)^{\frac{1}{2}}.
\end{aligned}
\label{4.50-4}
\end{equation}

Then by Definition \ref{def2.6-4}, Lemmas \ref{lem4.1-4}, \ref{lem4.3-4} and \ref{lem4.4-4}, it follows that $\psi_{t, T}\left(u^{1}, u^{2}\right)$ is a contractive function. For any constant $C_{4}>0$, taking $\tau=t-k+\frac{1}{2 \eta_{1}} \ln \frac{\|\tilde{u}_{\tau}\|_{C_{L^{2}(\Omega)}}^{2}+L\|\nabla \tilde{u}_{\tau}\|_{C_{L^{2}(\Omega)}}^{2}}{C_{4}}$, then it can be seen that the process $\{U(t, \tau)\}_{t \geq \tau}$ of problem (\ref{1.1-4}) is pullback $\mathcal{D}_{\eta_{1}}$-asymptotically compact in $C_{\mathcal H_{t}(\Omega)}$.$\hfill$$\Box$

From the above proofs, it follows the following theorem about the time-dependent pullback $\mathcal{D}_{\eta_{1}}$-attractors, which is one of the main results of this paper.

\begin{Theorem}\label{th4.1-2}
Under the assumptions of Lemma $\ref{lem4.3-2}$ and assume that the function $h(t)$ satisfies $(\ref{4.35-4})$, then there exists a unique time-dependent pullback $\mathcal{D}_{\eta_{1}}$-attractor $\mathcal A_{\eta_{1}}={}\left\{A_{\eta_{1}}(t): t \in \mathbb{R}\right\}$ of problem $(1.1)$ in $C_{\mathcal H_{t}(\Omega)}$.
\end{Theorem}
$\mathbf{Proof.}$ From Definitions \ref{def2.2-4} and \ref{def2.5-4}, Theorems $\ref{th3.1-4}-\ref{th3.2-4}$, Lemmas \ref{lem4.1-4}, \ref{lem4.3-4} and \ref{lem4.4-4}, it follows the existence and uniqueness of the above time-dependent pullback $\mathcal{D}_{\eta_{1}}$-attractor $\mathcal A_{\eta_{1}}$. $\hfill$$\Box$

\section{\large Regularity of the pullback attractors}
In this section, as in \cite{cp.4,qy.2,wzl.4,z.4}, we divide the weak solution $u$ of problem $(\ref{1.1-4})$ into two parts, and after finding the equations they satisfy respectively, we use the energy method to derive the regularity of $\mathcal A_{\eta_{1}}$ for the non-autonomous system $(\ref{1.1-4})$.
\begin{Theorem}
Under the assumptions of Theorem $\ref{th4.1-2}$, then $\mathcal A_{\eta_{1}}$ is bounded in $C_{\mathcal H_{t}^{1}(\Omega)}$.
\label{5.1-4}
\end{Theorem}

\noindent $\mathbf{Proof.}$ Since $L^{2}(\Omega) \subset H^{-1}(\Omega)$ is dense (see \cite{a.4,Evans}), for any $h(x, t) \in L_{\mathrm{loc}}^{2}(\mathbb R ; H^{-1}(\Omega))$, there exists a function ${h^\vartheta(x, t) } \in L_{\mathrm{loc}}^{2}(\mathbb R ; L^{2}(\Omega))$ such that
\begin{equation}
\|h-h^{\vartheta}\|<\vartheta,
\label{5.1-4}
\end{equation}
where $\vartheta > 0$ is a constant.

Fix $\tau \in \mathbb R$ and suppose $\phi \in \mathcal A_{\eta_{1}}$, then we can decompose the solution $U(t,\tau+\theta ){u(\tau+\theta) } = u(t)$ into the sum
$$
U(t, \tau+\theta) u(\tau+\theta)=U_{a}(t, \tau+\theta) u(\tau+\theta)+U_{b}(t, \tau+\theta) u(\tau+\theta),
$$
where $U_{a}(t, \tau+\theta) u(\tau+\theta)=v(t)$ and $U_{b}(t, \tau+\theta) u(\tau+\theta)=v_{1}(t)$ satisfy the following equations,  respectively,
\begin{equation}
\left\{\begin{array}{ll}
\partial_{t} v-\varepsilon(t)\partial_{t} \Delta v-a(l(u)) \Delta v=h(t)-h^{\vartheta}(t) & \text { in } \Omega \times(\tau, \infty), \\
v(x,t)=0 & \text { on } \partial \Omega \times(\tau, \infty), \\
v(x, \tau+\theta)=\phi(x,\theta), &\,\, x \in \Omega,
\label{5.2-4}
\end{array}\right.
\end{equation}
and
\begin{equation}
\left\{\begin{array}{ll}
\partial_{t} v_{1}-\varepsilon(t)\partial_{t} \Delta v_{1}-a(l(u)) \Delta v_{1}=f(u)+g(t,u_{t})+h^{\vartheta}(t) & \text { in } \Omega \times(\tau, \infty), \\
v_{1}(x,t)=0 & \text { on } \partial \Omega \times(\tau, \infty), \\
v_{1}(x, \tau+\theta)=0, &\,\, x \in \Omega.
\label{5.3-4}
\end{array}\right.
\end{equation}

Multiplying $(\ref{5.2-4})_{1}$ by $ - \Delta v$ and integrating it in $\Omega$, we arrive at
$$
\frac{d}{d t}\left(\|\nabla v\|^{2}+\varepsilon(t)\|\Delta v\|^{2}\right)+\left(2 a(l(u))-\varepsilon^{\prime}(t)\right)\|\Delta v\|^{2}=2(h-h^{\vartheta}(t),-\Delta v).
$$

By $(\ref{1.3-4})$, the Cauchy  and Poincar\'{e} inequalities, then put $t+\theta$ instead of $t$ in the obtained inequality, we deduce
\begin{equation}
\frac{d}{d t} I_{1}(t)+ \varrho_{1} I_{1}(t) \leq \vartheta^{2},
\label{5.4-4}
\end{equation}
where $I_{1}(t)=\|\nabla v\|_{C_{L^{2}(\Omega)}}^{2}+|\varepsilon_{t}|\|\Delta v\|_{C_{L^{2}(\Omega)}}^{2}$ and $0<\varrho_{1} \leq \frac{2+L}{\lambda_{1}^{-1}+|\varepsilon_{t}|}$.

Using the Gronwall inequality, we obtain
\begin{equation}
\|U_{a}(t, \tau+\theta) u(\tau+\theta)\|_{C_{\mathcal{H}_{t}^{1}(\Omega)}}^{2} \leq e^{-\varrho_{1}(t+k-\tau)}\left\|\phi\right\|_{C_{\mathcal{H}_{t}^{1}(\Omega)}}^{2}+\frac{\vartheta^{2}}{\varrho_{1}}.
\label{5.5-4}
\end{equation}

Similarly, multiplying $(\ref{5.3-4})_{1}$ by $ - \Delta v_{1}$ and integrating it in $\Omega$, we conclude
\begin{equation}
\frac{d}{d t}\left(\|\nabla v_{1}\|^{2}+\varepsilon(t)\|\Delta v_{1}\|^{2}\right)+\left(2 a(l(u))-\varepsilon^{\prime}(t)\right)\|\Delta v_{1}\|^{2}=2(f(u)+g(t,u_{t})+h^{\vartheta},-\Delta v_{1}).
\label{5.6-4}
\end{equation}

Besides, from $(\ref{1.6-4})$, assumptions (A1)$-$(A3) and the Young inequality, we derive
\begin{equation}
2(f(u),-\Delta v_{1}) \leq 2 \tilde{\eta}\|u\|^{2}+\frac{1}{2}\|\Delta v_{1}\|^{2},
\label{5.7-4}
\end{equation}
\begin{equation}
2(g(t, u_{t}),-\Delta v_{1}) \leq C_{g}\|u_{t}\|_{C_{L^{2}(\Omega)}}^{2}+\|\Delta v_{1}\|^{2}
\label{5.8-4}
\end{equation}
and
\begin{equation}
2(h^{\vartheta},-\Delta v_{1}) \leq\|h^{\vartheta}\|_{-1}^{2}+\frac{1}{2}\|\Delta v_{1}\|^{2}.
\label{5.9-4}
\end{equation}

Inserting (\ref{5.7-4})$-$(\ref{5.9-4}) into $(\ref{5.6-4})$ and using $(\ref{1.4-4})$, then putting $t+\theta$ instead of $t$ in the obtained equation, it can be seen from Lemma \ref{4.1-4} that there is $0<\varrho_{2}<\frac{1+L}{\lambda_{1}^{-1}+\left|\varepsilon_{t}\right|}$ such that
\begin{equation}
\frac{d}{d t} I_{2}(t)+\varrho_{2} I_{2}(t) \leq (2 \tilde{\eta}+C_{g}) R_{1}+2\|h^{\vartheta}\|_{-1}^{2},
\label{5.10-4}
\end{equation}
where $I_{2}(t)=\|\nabla v_{1}\|_{C_{L^{2}(\Omega)}}^{2}+|\varepsilon_{t}| \|\Delta v_{1}\|_{C_{L^{2}(\Omega)}}^{2}$ and $R_{1}=\rho^{2}(t)$.

Then by the Gronwall inequality, we conclude
\begin{equation}
\left\|U_{b}(t, \tau+\theta) u(\tau+\theta)\right\|^{2}_{C_{\mathcal H_{t}^{1}(\Omega)}}\leq R_{2},
\label{5.11-4}
\end{equation}
where $R_{2}=e^{-\varrho_{2} t} \int_{\tau}^{t} e^{\varrho_{2} s}[(2 \tilde{\eta}+C_{g}) R_{1}+2\|h^{\vartheta}\|_{-1}^{2}] d s.$.

Thanks to $(\ref{5.5-4})$ and $(\ref{5.11-4})$, for any $t \in \mathbb R$, we obtain
\begin{equation}
\operatorname{dist}\left(\mathcal{A}_{\eta_{1}}, \mathcal{\bar{B}}_{C_{\mathcal{H}_{t}^{1}(\Omega)}}\left(R_{2}\right)\right) \leq C e^{-\varrho(t+k-\tau)} \rightarrow 0, \,\,\,\, \tau \rightarrow-\infty,
\label{5.12-4}
\end{equation}
where ${\varrho>\varrho_{2}>0}$ and
\begin{equation}
\mathcal{\bar{B}}_{C_{\mathcal{H}_{t}^{1}(\Omega)}}\left(R_{2}\right)=\left\{u(t) \in \mathcal{\bar{B}}_{C_{\mathcal{H}_{t}^{1}(\Omega)}}:\|u(t)\|_{C_{\mathcal{H}_{t}^{1}(\Omega)}}^{2} \leq R_{2}\right\}.
\label{5.13-4}
\end{equation}
\qquad Consequently, we deduce that ${{\cal A}_{\eta_{1}}} \subseteq {\mathcal{\bar{B}}_{C_{\mathcal{H}_{t}^{1}(\Omega)}}\left(R_{2}\right)}$, which implies the time-dependent pullback $\mathcal D_{\eta_{1}}$-attractor ${{\cal A}_{\eta_{1}}}$ is bounded in $C_{\mathcal H_{t}^{1}(\Omega)}$. $\hfill$$\Box$

$\mathbf{Acknowledgment}$

This paper was partially supported by the NNSF of China with contract number 12171082 and by the fundamental funds for the central universities with contract number $2232022G$-$13$.

$\mathbf{Data\,\,availability\,\,statement}$

Data sharing not applicable to this article as no datasets were generated or analysed during the current study.

$\mathbf{Conflict\,\,of\,\,interest}$

The authors have no conflict of interest.

}
\end{document}